\documentclass[11pt]{amsart}
\usepackage[margin=28mm]{geometry}
\usepackage{amsmath,amssymb,amsthm,mathtools}
\usepackage[T1]{fontenc}
\usepackage{lmodern,microtype,enumitem}
\usepackage[hidelinks]{hyperref}
\newtheorem{theorem}{Theorem}[section]
\newtheorem{corollary}[theorem]{Corollary}
\newtheorem{proposition}[theorem]{Proposition}
\newtheorem{lemma}[theorem]{Lemma}

\newtheorem*{theoremA}{Theorem A}
\newtheorem*{theoremB}{Theorem B}
\newtheorem*{corollaryC}{Corollary C}
\newtheorem*{theoremD}{Theorem D}
\newtheorem*{conjectureE}{Conjecture E}
\theoremstyle{remark}
\newtheorem{remark}[theorem]{Remark}
\newcommand{\Area}{\operatorname{Area}}
\newcommand{\Ric}{\operatorname{Ric}}
\newcommand{\inj}{\operatorname{inj}}
\newcommand{\pr}{\operatorname{pr}}
\newcommand{\R}{\mathbb R}
\newcommand{\C}{\mathbb C}
\newcommand{\Z}{\mathbb Z}
\newcommand{\eps}{\varepsilon}
\newcommand{\one}{\mathbf 1}
\title{Sectional curvature and the rigidity of Perelman's width}
\author{Ziyang Qin}
\date{}
\begin{document}
\begin{abstract}
Perelman introduced the width $A(\alpha,g)$ of a nontrivial homotopy class
$\alpha$ of families of loops to prove finite extinction time of the Ricci flow
on certain three-manifolds. This paper deals with lower bounds and rigidity
results for this width under curvature bounds. First, we show that
$A(\alpha,g)\ge2\pi$ for every nontrivial class $\alpha$ on a closed manifold of
sectional curvature at most one, with equality for the class of degree one on
the round three-sphere. Second, we show that a metric on
the three-sphere with sectional curvature at most one and injectivity radius at
least $\pi$ is round if $A=2\pi$. This applies to quarter-pinched metrics. The
second result is reduced to a rigidity theorem for three sets in the round
three-sphere, whose proof uses a Gaussian double bubble theorem with unequal
interface weights.
\end{abstract}
\maketitle

\section{Introduction}\label{sec:intro}

Let $(M,g)$ be a closed Riemannian manifold, $\Lambda M$ the space of
contractible loops in $C^1(S^1,M)$, and $M\subset\Lambda M$ the constant
loops. For $c\in\Lambda M$, Perelman \cite{Perelman} defines $A(c,g)$ as the
infimum of the areas of the Lipschitz maps from the unit disk to $M$ spanning
$c$, and $A(\Gamma,g)=\sup_{c\in\Gamma}A(c,g)$ for $\Gamma\subset\Lambda M$.
For a nontrivial class $\alpha\in\pi_k(\Lambda M,M)$ let
\[
 A(\alpha,g)=\inf_{\Gamma\in\alpha}A(\Gamma,g).
\]
We call this number the \emph{$A$-width} of $\alpha$.

Along a Ricci flow on a closed three-manifold, the rate of change of
$A(\alpha,g(t))$ is at most $-2\pi-\frac12R_{\min}A(\alpha,g(t))$, where
$R_{\min}$ is the minimum of the scalar curvature \cite[Lemma~1.2]{Perelman}.
Together with the Ricci flow with surgery, this gives finite extinction time on
closed oriented three-manifolds whose prime decomposition has no aspherical
factors \cite{Perelman}. The constant term $-2\pi$ comes from the Gauss--Bonnet
theorem for a least area disk spanning a loop that moves by the curve shortening
flow.
Morgan and Tian gave a detailed account \cite[Chapters~18 and~19]{MorganTian}.

Perelman's proof uses a version of the least area disk argument of Hamilton
\cite[Section~11]{Hamilton} and a regularization of the curve shortening flow of
Altschuler and Grayson \cite{AltschulerGrayson}. Hamilton studied nonsingular
solutions of the normalized Ricci flow on closed three-manifolds. He showed that
the boundary tori of a persistent hyperbolic piece are incompressible
\cite[Sections~11 and~12]{Hamilton}. Otherwise some loop that is nontrivial on
such a torus would bound a disk outside the piece. Along the flow, for every $\delta>0$
the least area of such a disk would decrease at a rate of at least $2\pi-\delta$
for all large times. This is impossible, since areas are nonnegative. Here too the constant $2\pi$ comes from the Gauss--Bonnet
theorem.

In a parallel work \cite{CM}, Colding and Minicozzi gave another proof of
finite extinction time. They used a different notion of width, the width $W$ of
a nontrivial homotopy class of sweepouts by two-spheres. In \cite{CM08} they
showed that $W$ is the total area of finitely many branched minimal
two-spheres. The Gauss--Bonnet theorem on these spheres gives the constant term
$-4\pi$ in their estimate.

Zhou extended this approach to sweepouts by tori \cite{Zhou10} and by surfaces
of genus at least two \cite{Zhou17}. There the conformal structure of the domain
is an additional variable. Cheng and Zhou obtained sweepouts of closed connected
oriented four and five-manifolds by surfaces of genus at least two from
multisections \cite{ChengZhou}. They showed that every closed Riemannian four
or five-manifold contains a branched immersed closed minimal surface. In these
works the competitors are parametrized surfaces, and the minimal surfaces
obtained are branched immersions. The fillings in the definition of $A$ are
parametrized disks as well.

We write $K_g$ for the sectional curvature of $g$. For sweepouts, lower bounds
come from the critical points produced by min-max theory. If $K_g\le1$, a
branched minimal two-sphere has area at least $4\pi$, by the Gauss--Bonnet
theorem and the Gauss equation. Since $W$ is the total area of such spheres
\cite{CM08}, it follows that $W\ge4\pi$. Mazet and Rosenberg \cite{MR} studied
complete three-manifolds with $0\le K\le1$ that contain an embedded minimal
two-sphere of area $4\pi$. They proved that such a manifold is isometric to the
unit three-sphere or to a quotient of the Riemannian product of the unit
two-sphere with $\R$. Later, Mazet proved that the Simon--Smith width $\widetilde W(g)$
of a metric on $S^3$ with $K_g\le1$ is at least $4\pi$, with equality only for
the round metric \cite[Theorem~8]{Mazet}. Under $\Ric>0$ and scalar curvature at
least $6$, Marques and Neves proved $\widetilde W(g)\le4\pi$, with the same
equality case \cite[Theorem~1.1]{MN}.

We now restrict to $S^3$. We write $A(g)$ for the $A$-width $A(\alpha,g)$ of the
class $\alpha$ of degree one, and $W(g)$ for the width of \cite{CM08} for
sweepouts of degree one (Section~\ref{sec:widths}). Every circle on a two-sphere
bounds a disk of at most half the area of the sphere, hence $2A(g)\le W(g)$
(Lemma~\ref{lem:background}). The round metric $g_0$ of curvature one has
$A(g_0)=2\pi$ and $W(g_0)=4\pi$ (Corollary~\ref{cor:round}). A metric $g$ is
\emph{round of curvature $K$} if $Kg$ is isometric to $g_0$, and \emph{round} if
$K=1$.

The numbers $2\pi$ and $4\pi$ in the estimates of Perelman and of Colding and
Minicozzi are thus the widths of $g_0$. These estimates concern the change of
the widths along the flow. Among metrics on $S^3$ with $K_g\le1$, both $W$ and
$\widetilde W$ have minimum $4\pi$, and $\widetilde W$ attains it only for the
round metric. It is natural to ask whether $2\pi$ is the minimum of $A$, and
whether only the round metric attains it. Since $2A\le W$, the bound
for $A$ is the stronger one.

In the definition of the $A$-width, the functional maximized along a family of
loops is the least area of a filling of each loop. As in Perelman's argument,
our lower bounds are proved by the curve shortening flow and least area disks
spanning the loops. We do not explicitly use the existence of critical points of
a min-max problem. Throughout, classes in $\pi_k(\Lambda M,M)$ are taken without
base points (Section~\ref{sec:widths}).

Our first result is a lower bound for the $A$-width in every closed manifold.

\begin{theoremA}
Let $(M,g)$ be a closed Riemannian manifold with $K_g\le1$. Then
$A(\alpha,g)\ge2\pi$ for every nontrivial class
$\alpha\in\pi_k(\Lambda M,M)$, $k\ge1$. In particular every metric $g$ on
$S^3$ with $K_g\le1$ has
\[
 A(g)\ge2\pi,\qquad W(g)\ge4\pi .
\]
\end{theoremA}

By scaling, $K_g\le K$ with $K>0$ implies $A(\alpha,g)\ge2\pi/K$ for every
nontrivial $\alpha$. Since $2A\le W\le\widetilde W$ (Lemmas~\ref{lem:background}
and~\ref{lem:WMN}), Theorem~A also gives the lower bounds for $W$ and
$\widetilde W$ above.

Our second result is a rigidity theorem on $S^3$. Let $\inj(g)$ be the
injectivity radius of $g$.

\begin{theoremB}
Let $g$ be a metric on $S^3$ with $K_g\le1$ and $\inj(g)\ge\pi$. If
$A(g)=2\pi$, then $g$ is round.
\end{theoremB}

\begin{remark}\label{rem:inj}
If $K_g\le1$, then by Klingenberg's lemma $\inj(g)\ge\pi$ if and only if every
closed geodesic of $g$ has length at least $2\pi$ (Lemma~\ref{lem:kling}). So
the hypothesis of Theorem~B can also be stated as a condition on the lengths of
the closed geodesics.
\end{remark}

Theorem~B is the analogue for the $A$-width of the equality case of
\cite[Theorem~8]{Mazet}. Its last step uses Mazet's theorem that a complete
three-manifold with $K\le1$ containing an unstable totally geodesic immersed
unit two-sphere has constant curvature one \cite[Corollary~1]{Mazet}.

By the injectivity radius estimates of Klingenberg and Sakai
\cite{KlingenbergSakai} and of Cheeger and Gromoll \cite{CheegerGromoll},
Theorem~B applies to quarter-pinched metrics.

\begin{corollaryC}
Let $g$ be a metric on $S^3$ with $\frac14\le K_g\le1$. If $A(g)=2\pi$, then $g$
is round.
\end{corollaryC}

The proof of Theorem~B rests on a rigidity theorem for the round sphere, which we
now state. A \emph{$3$-fan} in $S^3$ consists of three distinct closed halves of
great two-spheres with a common boundary great circle. They are its \emph{pages},
and the three closed regions between them are its \emph{lunes}. For the round
metric, $|E|$ and $P(E)$ are the volume and the perimeter of $E\subset S^3$. Half
of $S^3$ has volume $\pi^2$, and a closed three-dimensional hemisphere and every
lune have perimeter $4\pi$. The following theorem describes three sets of
perimeter at most $4\pi$, each of less than half the volume of $S^3$, that cover
almost every point an odd number of times.

\begin{theoremD}
Let $E_1,E_2,E_3$ be sets of finite perimeter in $(S^3,g_0)$ with
$|E_i|<\pi^2$ and $P(E_i)\le4\pi$ for each $i$, such that
$\one_{E_1}+\one_{E_2}+\one_{E_3}$ is odd almost everywhere. Then
$E_1,E_2,E_3$ are, up to null sets, the three lunes of a $3$-fan in $S^3$.
\end{theoremD}

The name is taken from the theory of partitions of measures. B\'ar\'any and
Matou\v{s}ek \cite{BaranyMatousek} called a point of the plane together with $k$
rays from it a $k$-fan. They studied when several probability measures in the
plane can be divided by one $k$-fan in prescribed ratios. For their positive
results they transferred the problem to $S^2$, where a $k$-fan is given by $k$ great
semicircles from one point. Blagojevi\'c, Dimitrijevi\'c Blagojevi\'c and
McCleary \cite{BDM} defined a $k$-fan in $\R^d$ as $k$ distinct half-hyperplanes
whose common boundary is a linear subspace of dimension $d-2$, and a $k$-fan on
$S^{d-1}$ as its trace on the sphere. Up to taking closures and forgetting
orientation and labeling, their $3$-fans on $S^3$ are the $3$-fans above.

By Theorem~A, $A(g)\cdot\max K_g\ge2\pi$ for every metric on $S^3$
(Corollary~\ref{cor:round}). After normalizing $\max K_g=1$, the equality case
remains open if $\inj(g)<\pi$. We conjecture that the hypothesis on the
injectivity radius is not needed.

\begin{conjectureE}
Every metric $g$ on $S^3$ that is not round of any curvature satisfies
$A(g)\cdot\max K_g>2\pi$.
\end{conjectureE}

By Theorem~B, Conjecture~E is equivalent to the statement that $K_g\le1$ and
$\inj(g)<\pi$ imply $A(g)>2\pi$ (Proposition~\ref{prop:equiv}). The conjecture
holds for a class of doubly warped metrics and for some left-invariant metrics
(Corollary~\ref{cor:warped} and Theorem~\ref{thm:leftinv}). For smooth metrics
$g_n$ described by Mazet, $A(g_n)\cdot\max K_{g_n}\to2\pi$ while $g_n$ converges to
a nonround metric of class $C^{1,1}$ (Remark~\ref{rem:cylinder}).

To prove Theorem~B we compare $(S^3,g)$ with the round sphere. The map
$G=\exp_q\circ I\circ\exp_p^{-1}$, where $I$ is a linear isometry of tangent
spaces, is defined on the open ball of radius $\pi$ about $p$. By Rauch's theorem
it does not increase lengths. If $A(g)=2\pi$, modifications of $G$ that are
defined on all of $S^3$ carry almost optimal families of degree one to the round
sphere. For three nearby loops of such a family, the fillings of each pair and a
thin connecting annulus form a sphere. In the limit these spheres give three sets
of perimeter at most $4\pi$, each of at most half the volume, that cover almost
every point an odd number of times (Section~\ref{sec:round}). If one of them has
half the volume, it is a hemisphere, by Schmidt's inequality. Otherwise the three
sets are the lunes of a $3$-fan, by Theorem~D. Away from $-q$, the differential of
$G^{-1}$ preserves area on the boundary of the hemisphere, or on the pages of the
$3$-fan. The equality case of Rauch's theorem and the Codazzi equation then give a
totally geodesic unit sphere in $(S^3,g)$, after a change of the point $p$ if
necessary. If this sphere is unstable, $g$ is round by \cite[Corollary~1]{Mazet}.
Otherwise such stable spheres pass through every point of $S^3$, and
Proposition~\ref{prop:cover} excludes this.

For Theorem~D we first reduce to partitions of $S^3$ into three sets
(Lemma~\ref{lem:threecell}). We then apply a Gaussian double bubble theorem with
unequal interface weights (Theorem~\ref{thm:wgdb}) to the cones over these sets.
Theorem~\ref{thm:wgdb} is proved along the lines of the theorem of Milman and
Neeman \cite{MNg}. In $\R^n$, double bubbles with unequal weights were studied by
Lawlor \cite{Lawlor} (Remark~\ref{rem:lawlor}).

Section~\ref{sec:widths} collects preliminaries, and Section~\ref{sec:general}
proves Theorem~A. Section~\ref{sec:round} treats almost optimal families on the
round sphere, and Section~\ref{sec:threeset} proves Theorem~D. Theorem~B and
Corollary~C are proved in Section~\ref{sec:sharp}. Section~\ref{sec:conj}
concerns Conjecture~E.

\subsection*{Acknowledgments}
I would like to thank Xin Zhou for inspiring this work and for many
helpful discussions. I am also grateful to him for his careful reading of
earlier drafts, for many suggestions that greatly improved the paper, and
for his constant encouragement.

\subsection*{Use of AI tools}
I used generative AI tools to explore approaches to the problem and to edit the manuscript. I personally
checked every argument and take full responsibility for the entire paper.

\section{Preliminaries}\label{sec:widths}

This section fixes the notation, defines the widths of $S^3$ and compares them.

Throughout, $L(c)$ is the length of a loop. The area $\Area(f)$ of a Lipschitz or
$W^{1,2}$ map $f$ of a surface is the integral of
$|\partial_xf\wedge\partial_yf|\,dx\,dy$ in local coordinates. For a set
$E\subset S^3$, $|E|$ and $P(E)$ are the volume and the perimeter. We write
$\inj(p)$ for the injectivity radius at $p$ and $\inj(g)$ for its minimum. Closed
geodesics are nonconstant. The filling area $A(c,g)$ is defined as in the introduction also for
contractible Lipschitz loops $c$.

The set $\pi_k(\Lambda M,M)$ consists of the homotopy classes of maps of pairs
$(B^k,\partial B^k)\to(\Lambda M,M)$, without base points; no group structure is
used. A class is \emph{trivial} if it contains a map into $M$; then
$A(\alpha,g)=0$. A class of maps with base points lies in one of these classes, and is
trivial if that class is. Indeed, let $H$ be a homotopy of pairs from $\Gamma$ to
a map into $M$. Let $D_r\subset B^k\times[0,1]$, $0\le r\le1$, be disks with
boundary $\partial B^k\times\{0\}$ from $D_0=B^k\times\{0\}$ to
$D_1=\partial B^k\times[0,1]\cup B^k\times\{1\}$. The restrictions of $H$ to the
$D_r$ deform $\Gamma$ into $M$ relative to $\partial B^k$. So lower bounds for
$A(\alpha,g)$ also hold for classes with base points. Collapsing the circles over $\partial B^2$ turns $B^2\times S^1$
into $S^3$, and a map $\Gamma$ of pairs into the map $(x,t)\mapsto\Gamma(x)(t)$
of $S^3$ to $S^3$. This identifies $\pi_2(\Lambda S^3,S^3)$ with the homotopy
classes of maps $S^3\to S^3$, that is, with the degrees. Indeed, let
$F\colon B^2\times S^1\to S^3$ be continuous and constant on the circles over
$\partial B^2$, and $F_\delta(x,t)=\pr\int e_\delta(s)F(x,t-s)\,ds$ with
$\pr(y)=y/|y|$ and $e_\delta\ge0$ smooth, of integral one and supported in
$(-\delta,\delta)$. For small $\delta$, the map $x\mapsto F_\delta(x,\cdot)$ is continuous into
$C^1(S^1,S^3)$, and constant loops stay constant. The maps
$\pr((1-r)F+rF_\delta)$, $0\le r\le1$, join $F$ to $F_\delta$. If $F$ comes from a
map of pairs into $\Lambda S^3$, so do these maps. The same applies to homotopies
of such maps. Classes of nonzero degree are nontrivial, since a map into the
constant loops factors through $B^2$. A \emph{family of degree $d$} is a map of
pairs in the class of degree $d$. Reversing the loops changes the sign of the
degree and preserves filling areas, so $A(g)=A(\alpha,g)$ for either class
$\alpha$ of degree $\pm1$.

A \emph{sweepout} is a continuous map $\gamma\colon[0,1]\to C^0\cap
W^{1,2}(S^2,S^3)$ such that $\gamma(0)$ and $\gamma(1)$ are constant maps; it
induces a map $S^3\to S^3$, and $W(g)=\inf_\gamma\max_s\Area(\gamma(s))$ over
the sweepouts of degree $\pm1$. Sweepouts of the same degree are homotopic
through sweepouts. Indeed, take a homotopy of the induced maps $S^3\to S^3$. For small
$\delta$, mollify it at scale $\delta$ in the sphere variable as in the proof of
\cite[Lemma~D.1]{CM08}, and project it to $S^3$ by $\pr$. This joins the mollified
sweepouts $\gamma_\delta$. As $\gamma_\delta(s)\to\gamma(s)$ in $C^0\cap W^{1,2}$
uniformly in $s$, the maps $\pr((1-r)\gamma+r\gamma_\delta)$, $0\le r\le1$, join
$\gamma$ to $\gamma_\delta$. Reversing the parameter
changes the sign of the degree. So $W(g)$ is the width of the homotopy class of
sweepouts of degree one in the sense of \cite{CM08}. Colding and Minicozzi
define it with the energy in place of the area; the two widths agree
\cite[Proposition~1.2]{CM08}.

Filling areas are invariant under reparametrization. We have
$A(c\circ m,g)=A(c,g)$ for every Lipschitz weakly monotone map $m$ of the circle
of degree one. Indeed, for a lift of $m$, the annulus
$(t,r)\mapsto c((1-r)m(t)+rt)$ has zero area. A reflection of the disk reverses
the orientation of the boundary. Joining
two loops at uniform distance $d\le\frac12\inj(g)$ by minimizing geodesics
gives, with $C$ depending only on $g$,
\begin{equation}\label{eq:Acont}
 |A(c,g)-A(c',g)|\le C\,d\,\bigl(L(c)+L(c')\bigr).
\end{equation}
In particular $A(\cdot,g)$ is continuous on $\Lambda M$.

The next lemma collects the elementary properties of the widths and compares
them.

\begin{lemma}\label{lem:background}
For metrics on $S^3$ the widths are finite, and
\[
 2A(g)\le W(g),\qquad A(g_0)\le2\pi,\qquad W(g_0)\le4\pi .
\]
They are monotone in the metric tensor, and $A(tg)=tA(g)$, $W(tg)=tW(g)$ for
$t>0$.
\end{lemma}
\begin{proof}
Monotonicity and homogeneity follow by comparing area integrands. Let $\gamma$
be a sweepout of degree $\pm1$ and $\eps>0$. At the cost of $\eps$ in $\max_s\Area(\gamma(s))$, the sweepout $\gamma$ may be
assumed continuous into $C^2(S^2,S^3)$, with $\gamma(0)$ and $\gamma(1)$ still
constant maps. This follows from \cite[Lemma~D.1]{CM08} and the continuity of the
area in $W^{1,2}$. Let $\ell_a(t)=(\sin a\,e^{it},\cos a)$, $0\le a\le\pi$, be the latitude
circles of $S^2\subset\C\times\R$. The loops $\Gamma(s,a)=\gamma(s)\circ\ell_a$
form a family parametrized by $[0,1]\times[0,\pi]$, with a constant loop at
each boundary parameter. Its map $S^3\to S^3$ is the map of $\gamma$ composed with a
map of degree one, so it has degree $\pm1$. The loop $\Gamma(s,a)$
bounds the images under $\gamma(s)$ of the two caps of $S^2$ bounded by
$\ell_a$, whose areas add up to $\Area(\gamma(s))$. Hence
$A(\Gamma(s,a),g)\le\frac12\Area(\gamma(s))$ and $2A(g)\le W(g)+\eps$. The
sweepout of $S^3$ by the spheres $\{x_4=t\}$ gives $W(g)<\infty$ and
$W(g_0)\le4\pi$, hence $A(g_0)\le2\pi$.
\end{proof}

The Simon--Smith width $\widetilde W(g)$ is the infimum, over smooth families
$F_t$, $-1\le t\le1$, of diffeomorphisms of $S^3$ isotopic to the identity, of
$\max_t\Area(F_t(\{x_4=t\}))$ \cite{MN,Mazet}. The next lemma shows that it is
at least $W(g)$.

\begin{lemma}\label{lem:WMN}
Every metric $g$ on $S^3$ satisfies $W(g)\le\widetilde W(g)$.
\end{lemma}
\begin{proof}
Parametrizing the spheres $\{x_4=t\}$ by $S^2$, each family $F_t(\{x_4=t\})$ is
a sweepout of degree one.
\end{proof}

\section{The lower bound}\label{sec:general}

In this section we prove Theorem~A. Throughout, $(M,g)$ is a closed Riemannian
manifold. We follow Perelman's argument \cite{Perelman} for a fixed metric. The
loops of a family are first smoothed and lifted to ramps; compare
\cite[3.2]{Perelman}. Ramps are curves in the product of $M$ with a circle
$S^1_\lambda$ of small length $\lambda$. They wind once around the circle and are
graphs over it. The precise definition is given below. The curve shortening flow
of ramps exists for all time, and its projection to $M$ deforms the loops. The
filling areas are compared with the solution of a linear differential equation.
The result is the following proposition.

\begin{proposition}\label{prop:decay}
Let $(M,g)$ be a closed Riemannian manifold with $K_g\le1$, let $\Gamma\colon P\to\Lambda M$ be a continuous map from a
compact metric space, and $T,\xi>0$. There is a continuous deformation $\Gamma^t$,
$-1\le t\le T$, of $\Gamma$, with $\Gamma^t=\Gamma$ for $t=-1$, through
continuous maps $P\to\Lambda M$ such that
\begin{enumerate}[label=\textup{(\alph*)}]
\item if $\Gamma(x)$ is a constant loop, then $\Gamma^t(x)=\Gamma(x)$ for all
$t$;
\item for every $x\in P$ either $L(\Gamma^T(x))\le\xi$, or
\[
 A(\Gamma^T(x),g)\le w(T)+\xi,\qquad\text{where }w'=w-2\pi,\ w(0)=A(\Gamma(x),g).
\]
\end{enumerate}
\end{proposition}

Proposition~\ref{prop:decay} is the analogue of \cite[3.1]{Perelman} for a fixed
metric. There the metric is a Ricci flow on a three-manifold, and
$w'=-2\pi-\frac12R_{\min}w$. Here $w(T)=2\pi+e^T(w(0)-2\pi)$. Before the proof we
show that the proposition implies Theorem~A.

\begin{proof}[Proof of Theorem~A]
Let $\Gamma$ represent $\alpha$ with $a_0=A(\Gamma,g)<2\pi$. Choose $T$ with
$2\pi+e^T(a_0-2\pi)<0$, and $\xi>0$ smaller than the absolute value of this
number and than $\inj(g)$. Filling areas are nonnegative, so every loop
$c_x=\Gamma^T(x)$ has length at most $\xi$. By (a), $\Gamma^t$ is a homotopy of
maps of pairs, and
\[
 (x,s)\longmapsto\Bigl(y\mapsto
 \exp_{c_x(0)}\bigl((1-s)\exp_{c_x(0)}^{-1}c_x(y)\bigr)\Bigr),\qquad0\le s\le1,
\]
deforms $\Gamma^T$ into the constant loops. Hence $\alpha$ is trivial. The
statement for $S^3$ follows with a class of degree one and $2A\le W$
(Lemma~\ref{lem:background}).
\end{proof}

Theorem~A and Lemma~\ref{lem:background} determine the widths of the round metric.
Theorem~A also gives a bound for every metric on $S^3$.

\begin{corollary}\label{cor:round}
The round metric has $A(g_0)=2\pi$ and $W(g_0)=4\pi$. Every metric $g$ on $S^3$
has $\max K_g>0$ and $A(g)\cdot\max K_g\ge2\pi$.
\end{corollary}
\begin{proof}
Lemma~\ref{lem:background} gives the upper bounds and Theorem~A the lower ones.
If $K_g\le0$, Theorem~A gives $A(g)\ge2\pi/K$ for every $K>0$, which is
impossible, as $A(g)<\infty$. The last inequality is Theorem~A for the metric
$(\max K_g)\,g$.
\end{proof}

The rest of this section proves Proposition~\ref{prop:decay}, for a small fixed
$\lambda$. The filling area $F$ of the ramps in $M\times\R^2$ satisfies
$F'\le F-2\pi$ up to an error that is small unless the ramp is steep
(Lemma~\ref{lem:diff}). The times at which the ramp is both long and steep have
small measure (Lemma~\ref{lem:steep}).

We begin with the definition of ramps. Let $S^1_\lambda\subset\R^2$ be the circle
of length $\lambda$ about the origin, and
\[
 M_\lambda=M\times S^1_\lambda\subset M\times\R^2
\]
with the product metrics. Let $U$ be the unit vector field on $M_\lambda$
tangent to the circles, $n$ the outward unit normal of $M_\lambda$ in
$M\times\R^2$, and $p\colon M\times\R^2\to M$ the projection. Then $U$ is
parallel on $M_\lambda$, and $\nabla_UU=-\frac{2\pi}\lambda n$ in
$M\times\R^2$. For orthonormal $X,Y$ with components $X_1,Y_1$ tangent to $M$,
the curvature of the plane $X\wedge Y$ is $\langle R^g(X_1,Y_1)Y_1,X_1\rangle$.
This is zero if $X_1,Y_1$ are linearly dependent. So the sectional curvatures of $M\times\R^2$ and of $M_\lambda$ satisfy
\begin{equation}\label{eq:secprod}
 K\le K_+:=\max\bigl(\{0\}\cup\{K_g(P)\}\bigr),
\end{equation}
with $P$ ranging over the tangent planes of $M$. In particular $K\le1$ if
$K_g\le1$.

A \emph{ramp} \cite[2.3]{Perelman} is a smooth closed curve $c$ in
$M_\lambda$ that winds once around the circle factor and whose unit tangent $S$
has $u=\langle S,U\rangle>0$. A ramp is embedded, it is the graph of a map
$S^1_\lambda\to M$, and
\begin{equation}\label{eq:rise}
 \int_cu\,ds=\lambda .
\end{equation}
Let $\vec H$ and $\vec H'$ be the curvature vectors of $c$ in $M_\lambda$ and in
$M\times\R^2$, and $k=|\vec H|$. The second fundamental form of $M_\lambda$ on $S$
is $u^2\nabla_UU$, so
\begin{equation}\label{eq:kappa}
 \vec H=\vec H'+\frac{2\pi}\lambda\,u^2n .
\end{equation}
For a Lipschitz loop $c$ in $M$, parametrized by $x\in\R/\Z$, let
$c_\lambda(x)=(c(x),\lambda x)$, where $\lambda x$ is the point of
$S^1_\lambda$ at arclength $\lambda x$. Its length is at most $L(c)+\lambda$. If
$c$ is smooth, $c_\lambda$ is a ramp, the \emph{ramp over $c$}. The ramp over a
constant loop is a circle, a geodesic of $M_\lambda$.

The next lemma shows that the curve shortening flow of a ramp exists for all time
and consists of ramps, with bounds on its curvature.

\begin{lemma}[the flow of ramps]\label{lem:ramp}
Let $(M,g)$ be a closed Riemannian manifold. The curve shortening flow $c^t$ in $M_\lambda$ of a ramp $c^0$ exists for all
$t\ge0$, and every $c^t$ is a ramp. With $s$ the arclength,
\begin{enumerate}[label=\textup{(\alph*)}]
\item the length satisfies $\frac d{dt}L(c^t)=-\int_{c^t}k^2\,ds$, and
$|\partial_su|\le k$;
\item the minimum of $u$ is nondecreasing in $t$, and
$\max_{c^t}(k/u)\le e^{K_+t}\max_{c^0}(k/u)$, with $K_+$ as in \eqref{eq:secprod};
\item as graphs $c^t=\{(y(z,t),z)\}$ over $z\in S^1_\lambda$, the ramps solve
$\partial_ty=\nabla_zy_z/(1+|y_z|^2)$, and on a set of initial ramps that is
compact in $C^2$ the solution depends continuously on $y(\cdot,0)$ in
$C^0([0,T],C^1)$ for every $T$.
\end{enumerate}
\end{lemma}
\begin{proof}
The argument is that of \cite[\S2]{AltschulerGrayson} and \cite[2.3]{Perelman}.
As long as the flow exists, $[\partial_t,\partial_s]=k^2\partial_s$, and since $U$ is
parallel,
\[
 \partial_tu=\partial_s^2u+k^2u,\qquad
 \partial_tk^2=\partial_s^2k^2-2|\nabla_s^\perp \vec H|^2+2k^4
 +2\langle R(\vec H,S)S,\vec H\rangle .
\]
The maximum principle gives $u>0$ with $\min u$ nondecreasing. The identity
$\partial_su=\langle \vec H,U\rangle$ gives $|\partial_su|\le k$, and the first
variation of length gives (a). At a positive maximum of $v=k^2/u^2$ one has
$\partial_sk^2=2k^2\partial_su/u$, and
\[
 \partial_tv=\partial_s^2v
 +\frac2{u^2}\Bigl(\frac{k^2(\partial_su)^2}{u^2}-|\nabla_s^\perp \vec H|^2\Bigr)
 +\frac2{u^2}\langle R(\vec H,S)S,\vec H\rangle .
\]
The bracket is nonpositive there, since $|\partial_sk^2|\le2k|\nabla_s^\perp \vec H|$.
The last term is at most $2K_+v$ by \eqref{eq:secprod}. So $\max v\le
e^{2K_+t}\max v|_{t=0}$, which is (b). Thus $k$ stays bounded on bounded time
intervals, and the flow exists for all time \cite[Propositions~3.2
and~3.11]{Smoczyk}. For curves in $\R^3$ see
\cite[Theorem~1.13]{AltschulerGrayson}. The vector field
$\nabla_zc_z/|c_z|^2$ along the graph has normal part $\vec H$ and no component
along the circle, which gives the equation in (c). There
$u^{-2}=1+|y_z|^2$ and $|\nabla_zy_z|\le k/u^3$, so (b) bounds $y$ in $C^2$
uniformly on the given set, and the equation is uniformly parabolic. For two
solutions, in an isometric embedding $M\subset\R^m$, the difference
$d=y_1-y_2$ satisfies $\frac d{dt}\int|d|^2\le C\int|d|^2$, and continuity in
$C^1$ follows by interpolation.
\end{proof}

The error term in Lemma~\ref{lem:diff} below is at most $2\pi\max u$, so we have
to bound the times at which a long ramp is steep. The next lemma shows that such
a ramp loses length at a definite rate.

\begin{lemma}\label{lem:steep}
Let $c$ be a ramp and $0<\delta<1$. If $\max_cu>\delta$ and
$L(c)>2\lambda/\delta$, then $\int_ck^2\,ds\ge\delta^3/(8\lambda)$.
\end{lemma}
\begin{proof}
By \eqref{eq:rise}, $\min u\le\lambda/L(c)<\delta/2$, so there is an arc $J$
on which $\delta/2\le u\le\delta$, with $u=\delta/2$ at one end and $u=\delta$
at the other, and $|J|\,\delta/2\le\lambda$ by \eqref{eq:rise}. By
Lemma~\ref{lem:ramp}(a) and the Cauchy--Schwarz inequality,
$\delta/2\le\int_Jk\,ds\le|J|^{1/2}(\int_Jk^2)^{1/2}$.
\end{proof}

Since $\frac d{dt}L=-\int k^2\,ds$, the lemma shows that a flow of ramps longer
than $2\lambda/\delta$ has $\max u>\delta$ only on a set of times of measure at
most $8\lambda L(c^0)/\delta^3$.

We now compare filling areas in $M$ and in $M\times\R^2$. For a Lipschitz curve
$c$ in $M\times\R^2$ such that $p\circ c$ is contractible in $M$, let $F(c)$ be
its filling area in $M\times\R^2$. It is defined as $A(c,g)$ with $M\times\R^2$
in place of $M$. The ramps are filled in $M\times\R^2$ instead of being projected to $M$.
There a ramp is embedded for every $\lambda$. If its projection to $M$ is
contractible, part (c) of the next lemma gives a least area disk of class $C^2$
up to the boundary, to which Lemma~\ref{lem:GB} applies if $K_g\le1$. Part (a) transfers bounds on $F$ to the projected loops, and part (b)
bounds $F$ at the initial time.

\begin{lemma}\label{lem:F}
Let $(M,g)$ be a closed Riemannian manifold, $c$ a ramp in $M_\lambda$ and $b$ a loop in $M$.
\begin{enumerate}[label=\textup{(\alph*)}]
\item The projection satisfies $L(p\circ c)\le L(c)$, and $A(p\circ c,g)\le
F(c)$ if $p\circ c$ is contractible.
\item If $b$ is Lipschitz and contractible, then $F(b_\lambda)\le
A(b,g)+\frac\lambda{2\pi}L(b)+\frac{\lambda^2}{4\pi}$.
\item If $p\circ c$ is contractible, there is a nonconstant map
$f\in C^2(\overline D,M\times\R^2)$, conformal and harmonic in $D$, whose
restriction to $\partial D$ is a weakly monotone parametrization of $c$ of
degree one, with $\Area(f)=F(c)$.
\end{enumerate}
\end{lemma}
\begin{proof}
(a) The projection does not increase lengths or areas. (b) Add to a disk in
$M\times\{0\}$ spanning $b$ the annulus $(x,r)\mapsto(b(x),re^{2\pi ix})$,
$0\le r\le\lambda/2\pi$. Its coordinate fields are orthogonal, of lengths at
most $|b'|+2\pi r$ and $1$. (c) The manifold $M\times\R^2$ is homogeneously
regular, and $c$ is a smooth Jordan curve, contractible in $M\times\R^2$. By
Morrey's solution of the problem of Plateau \cite{Morrey}, it bounds a disk $f$,
conformal and harmonic, of least area. Here the competitors are the continuous
maps of class $W^{1,2}$ whose boundary values are weakly monotone
parametrizations of $c$. By
the boundary regularity of Heinz and Hildebrandt
\cite[Theorem, (iii)]{HeinzHildebrandt}, $f$ is of class $C^2$ up to the boundary
(cf.\ \cite[proof of Lemma~19.2]{MorganTian}). Lipschitz disks spanning $c$ are
competitors, and $f$ is one of them up to a reparametrization of the boundary
as in Section~\ref{sec:widths}. So $\Area(f)=F(c)$.
\end{proof}

The next lemma is the boundary inequality that follows from the Gauss--Bonnet
theorem and the Gauss equation, for conformal harmonic disks. It allows branch
points.

\begin{lemma}\label{lem:GB}
Let $N$ be a Riemannian manifold with $K\le1$, $c$ a smooth embedded closed
curve in $N$ with curvature vector $\vec H'$, and $f\in C^2(\overline D,N)$
nonconstant, conformal and harmonic in $D$, with $f(\partial D)\subset c$. Then
\[
 \int_0^{2\pi}\bigl\langle \vec H'(f(e^{it})),\,\partial_rf(e^{it})
 \bigr\rangle\,dt\ \le\ -2\pi+\Area(f),
\]
with polar coordinates $re^{it}$ on $D$.
\end{lemma}
\begin{proof}
Let $q=|f_x|^2=|f_y|^2$, so that $\Area(f)=\int_Dq$. A harmonic map that is
constant on an open set is constant \cite{Sampson}, so $\{q>0\}$ is dense in
$D$. There $f$ is a minimal immersion with metric $q\,|dw|^2$ and curvature
$-\frac1{2q}\Delta\log q\le1$, that is, $\Delta q\ge|\nabla q|^2/q-2q^2$. For
$b>0$ this gives $\Delta\log(q+b)\ge-2q$ on $\{q>0\}$, hence on $D$.
Integrating over disks of radius less than one and passing to the limit,
\begin{equation}\label{eq:GB1}
 \int_0^{2\pi}\frac{\partial_rq}{2(q+b)}\,dt\ \ge\ -\Area(f).
\end{equation}
The relations $|f_r|^2=|f_t|^2=q$, $\langle f_r,f_t\rangle=0$ and
$\nabla_rf_r+f_r+\nabla_tf_t=0$ hold on $\partial D$. Let $s(t)$ be the arclength
coordinate of $f(e^{it})$ on $c$ and $S$ the unit tangent of $c$. Then
$f_t=s'S$, $\nabla_tf_t=s''S+s'^2\vec H'$, and
$\langle S,f_r\rangle=0$ (where $s'=0$ also $f_r=0$). Hence
$\frac12\partial_rq=-q(1+\langle \vec H'\circ f,f_r\rangle)$ on $\partial D$.
Insert this in \eqref{eq:GB1} and let $b\to0$. By dominated convergence,
\[
 \int_{\{t\,:\,q(e^{it})>0\}}\bigl(1+\langle \vec H'\circ f,f_r\rangle
 \bigr)\,dt\ \le\ \Area(f).
\]
By Heinz and Hildebrandt \cite[Theorem, (ii)]{HeinzHildebrandt}, at a zero
$w_0\in\partial D$ of $df$ one has $\partial_wf=a(w-w_0)^l+o(|w-w_0|^l)$ in
local coordinates, with an integer $l\ge1$ and $a\ne0$. So the zeros of $q$ on
$\partial D$ are isolated, hence finite in number.
\end{proof}

\begin{remark}
For an immersed disk, Lemma~\ref{lem:GB} is the Gauss--Bonnet theorem with the
Gauss equation $K_f\le1$ for minimal surfaces. For an embedded closed geodesic it
gives $\Area(f)\ge2\pi$. Branch points are needed only in
Lemma~\ref{lem:barrier}. The interior zeros of $q$ are handled by $\log(q+b)$, and
those on the boundary by \cite{HeinzHildebrandt}. The disks of Lemma~\ref{lem:F}(c) are in fact immersions. Their
component $v$ in $\R^2$ is harmonic, with $|v|=\lambda/2\pi$ on $\partial D$ and
$|v|<\lambda/2\pi$ in $D$, so $\partial_r|v|^2>0$ on $\partial D$ by Hopf's
lemma, and $f_r\ne0$ there. By conformality and as $u>0$, $v$ maps $\partial D$
diffeomorphically onto a circle. By the theorems of Rad\'o, Kneser and Choquet
and of Lewy \cite[Sections~3.1 and~2.2]{Duren}, $v$ is a diffeomorphism of $D$
onto a disk, so $df\ne0$ in $D$.
\end{remark}

The next lemma gives the differential inequality for the filling area of a flow
of ramps. It is the computation of \cite[1.2]{Perelman} for ramps in a fixed
metric. Without the last term it is the differential inequality $F'\le F-2\pi$,
with the comparison equation $w'=w-2\pi$ of Proposition~\ref{prop:decay}. The
last term comes from the curvature $2\pi/\lambda$ of $S^1_\lambda$ in $\R^2$,
through \eqref{eq:kappa}.

\begin{lemma}\label{lem:diff}
Let $(M,g)$ be a closed Riemannian manifold with $K_g\le1$, let $c^t$ be a flow
of ramps with $p\circ c^0$ contractible, and $F(t)=F(c^t)$. Then $F$ is locally Lipschitz, and
\[
 \limsup_{h\downarrow0}\frac{F(t+h)-F(t)}h\ \le\
 F(t)-2\pi+\frac{2\pi}\lambda\int_{c^t}u^2\,ds
 \ \le\ F(t)-2\pi\bigl(1-\max_{c^t}u\bigr).
\]
\end{lemma}
\begin{proof}
The strip swept by the ramps between two times has area at most a constant
times the difference of the times, so $F$ is locally Lipschitz. Fix $t$, and let
$f$ be the disk of Lemma~\ref{lem:F}(c) for $c^t$. Let $X$ be a smooth vector
field with compact support on $M\times\R^2$ that equals $\vec H$ along $c^t$, with
flow $e^{hX}$. Since $f$ is conformal and harmonic and $C^1$ up to the boundary,
\[
 \frac d{dh}\Big|_{h=0}\Area(e^{hX}\circ f)=\int_{\partial D}\langle X\circ f,
 f_r\rangle .
\]
The integral is with respect to arclength on $\partial D$. The formula holds also
at branch points, since the area integrand is $q$ times a function smooth in $h$.
Parametrize the ramps so that the velocity is $\vec H$, which does not change
$F$. Then the curves $e^{hX}\circ c^t$ and $c^{t+h}$ are at distance $O(h^2)$, so
$F(t+h)\le\Area(e^{hX}\circ f)+O(h^2)$. Hence the left side of the lemma is at
most $\int_{\partial D}\langle\vec H\circ f,f_r\rangle$. By
\eqref{eq:kappa}, Lemma~\ref{lem:GB} and \eqref{eq:secprod},
\[
 \int_{\partial D}\langle\vec H\circ f,f_r\rangle
 \le-2\pi+F(t)+\frac{2\pi}\lambda\int_{\partial D}u^2|f_r| .
\]
Finally, on $\partial D$, $|f_r|$ is the length of the tangential derivative of
$f$. So $|f_r|$ times arclength pushes forward to $ds$ under the monotone map
$f|_{\partial D}$. Moreover $\int u^2ds\le\lambda\max u$ by \eqref{eq:rise}.
\end{proof}

We now combine these lemmas to prove Proposition~\ref{prop:decay}.

\begin{proof}[Proof of Proposition~\ref{prop:decay}]
We first smooth the loops. Embed $M$ isometrically in $\R^m$ with nearest point
projection $\pr$, and let $\Gamma_h(x)=\pr(\Gamma(x)*m_h)$ for a mollifier $m_h$
of width $h$ on the circle. For small $h$ this is a homotopy of maps
$P\to\Lambda M$ fixing constant loops, and for $h>0$ it is continuous in $x$
with values in $C^2$. Fix by \eqref{eq:Acont} an $h>0$ with
$A(\Gamma_h(x))\le A(\Gamma(x))+\xi_0$ for all $x$, where $\xi_0=e^{-T}\xi/4$. Let $L_0$ bound the lengths of
the loops $\Gamma_h(x)$. For $-1\le t\le0$ let $\Gamma^t=\Gamma_{(1+t)h}$, with
$\Gamma_0=\Gamma$.

We next deform the loops. Let $c_x^t$ be the flow of ramps in $M_\lambda$
starting from the ramp over $\Gamma_h(x)$, and $\Gamma^t(x)=p\circ c_x^t$ for
$0\le t\le T$, parametrized by the circle factor. By Lemma~\ref{lem:ramp}(c)
this is continuous in $(x,t)$ with values in $C^1$, and constant loops stay
fixed.

It remains to estimate the filling areas. Fix $x$ and $\delta\in(0,1)$, and write
$c^t=c^t_x$, $F(t)=F(c^t)$ and $a=A(\Gamma(x))$. By Lemma~\ref{lem:F}(b),
\begin{equation}\label{eq:F0}
 F(0)\le a+\xi_0+\tfrac\lambda{2\pi}L_0+\tfrac{\lambda^2}{4\pi},\qquad
 L(c^0)\le L_0+\lambda .
\end{equation}
The length is nonincreasing. If $L(c^T)\le2\lambda/\delta$, then
$L(\Gamma^T(x))\le2\lambda/\delta$. Otherwise $L(c^t)>2\lambda/\delta$ on
$[0,T]$, and the set $\mathcal B$ of times with $\max_{c^t}u>\delta$ has, by
Lemmas~\ref{lem:steep} and~\ref{lem:ramp}(a),
\[
 |\mathcal B|\,\frac{\delta^3}{8\lambda}\le\int_0^T\!\!\int_{c^t}k^2\,ds\,dt
 \le L(c^0),
 \qquad\text{so}\qquad
 |\mathcal B|\le\frac{8\lambda(L_0+\lambda)}{\delta^3}.
\]
By Lemma~\ref{lem:diff} and comparison with the linear equation,
\[
 F(T)\le e^TF(0)-2\pi\int_0^Te^{T-t}\bigl(1-\max_{c^t}u\bigr)\,dt
 \le 2\pi+e^T\bigl(F(0)-2\pi\bigr)+2\pi e^T\bigl(\delta+|\mathcal B|\bigr).
\]
With \eqref{eq:F0} and Lemma~\ref{lem:F}(a),
\[
 A(\Gamma^T(x))\le F(T)\le w(T)
 +e^T\Bigl(\xi_0+\tfrac\lambda{2\pi}L_0+\tfrac{\lambda^2}{4\pi}+2\pi\delta
 +\frac{16\pi\lambda(L_0+\lambda)}{\delta^3}\Bigr).
\]

Choose $\delta$ with $2\pi\delta e^T\le\xi/4$. Then choose $\lambda$ so small
that $2\lambda/\delta\le\xi$ and the terms containing $\lambda$ add up to at
most $e^{-T}\xi/2$. The constants do not depend on $x$.
\end{proof}

\begin{remark}\label{rem:perelman}
Perelman \cite[\S3]{Perelman} and Morgan--Tian \cite[Chapter~19]{MorganTian}
project the ramps to $M$ and let $\lambda\to0$. Along a Ricci flow the term
$-\frac12R_{\min}A$ in Perelman's estimate \cite[Lemma~1.2]{Perelman} comes from
the identity
$\operatorname{tr}_P\Ric=K(P)+\frac12R$ for planes $P$ tangent to a
three-manifold. For planes in $M\times\R^2$ it would need an additional error
estimate. In a fixed metric \eqref{eq:secprod} replaces it.
\end{remark}

\section{Almost optimal families on the round sphere}\label{sec:round}

In this section we study families of degree one on the round sphere whose filling
areas are close to $2\pi$. We take three loops of such a family that are close to
each other and glue their fillings in pairs, so that three spheres arise. Each
sphere is described by an integer-valued function of bounded variation
(Lemma~\ref{lem:bounding}). The set where this function is odd has perimeter at
most the area of the sphere. Theorem~\ref{thm:switching} selects three
nearby loops, and for each pair of them an odd set of the corresponding sphere.
These three sets cover each point of $S^3$ an odd number of times, outside a set of
small volume.
Proposition~\ref{prop:limits} passes to the limit, and
Proposition~\ref{prop:lune} describes the case in which a limit is a lune.

In this section and in Section~\ref{sec:threeset} the sphere $S^3\subset\R^4=\C^2$
carries the round metric $g_0$, its volume form $\mathrm{vol}$ and an orientation
fixed below. The notations $|E|$, $P(E)$, $\mathcal H^k$, $\operatorname{Jac}$ and
$\Area$ refer to $g_0$. For $u\in L^1(S^3)$ and an open set $W$, $|Du|(W)$ is the supremum of
$\int_Wu\operatorname{div}X$ over $C^1$ vector fields $X$ with compact support in
$W$ and $|X|\le1$. If $|Du|(S^3)<\infty$, then $u\in BV(S^3)$, and $Du=n_u|Du|$,
$|n_u|=1$, is the vector measure with $\int u\operatorname{div}X=-\int\langle
X,dDu\rangle$. We put $P(E)=|D\one_E|(S^3)$, and $BV(S^3;\Z)$ consists of the
integer-valued elements of $BV(S^3)$. The standard facts on $BV$ functions in
$\R^3$ \cite{AFP,Maggi} carry over through charts. Sequences bounded in $BV(S^3)$
have subsequences converging in $L^1$ and almost everywhere. If $u_n\to u$ in
$L^1$, then $|Du|(W)\le\liminf|Du_n|(W)$ for open $W$. If moreover
$\sup_n|Du_n|(S^3)<\infty$, then $Du_n\to Du$ weakly-*
\cite[Theorem~3.23, Proposition~3.13]{AFP}.
A function $u$ with $|Du|(W)=0$ on a connected open set $W$ is almost everywhere
constant on $W$. For $u\in BV(S^3;\Z)$ the sets $\{u\ge k\}$, $k\in\Z$, have finite
perimeter and $|Du|=\sum_k|D\one_{\{u\ge k\}}|$ as measures (coarea formula
\cite[Theorem~3.40]{AFP}). Currents are those of \cite{Federer} in $\R^4$, and for
$u\in L^1(S^3)$ the $3$-current $u[\![S^3]\!]$ integrates $3$-forms over $S^3$
with the density $u$.

We use Schmidt's isoperimetric inequality on spheres with its equality case
\cite{Schmidt}, \cite[(2.2)]{BDF}. A set $E\subset S^3$ of finite perimeter with
$|E|=\pi(2r-\sin2r)$, $0<r<\pi$, has $P(E)\ge4\pi\sin^2r$, with equality only for
geodesic balls of radius $r$, up to null sets. In particular a set of volume
$\pi^2$ and perimeter at most $4\pi$ is a closed hemisphere up to a null set.
Since $\pi(2r-\sin2r)\le\frac{4\pi}3r^3$ and $4\pi\sin^2r\ge16r^2/\pi$ for
$r\le\frac\pi2$,
\begin{equation}\label{eq:schmidt}
 |E|\le c\,P(E)^{3/2}\qquad\text{if }|E|\le\pi^2,\qquad c=\pi^{5/2}/48 .
\end{equation}

The \emph{pages} of a great circle $C$ are the closed halves of great two-spheres
bounded by $C$. If $C=\{z_2=0\}$, they are $\{(z_1,re^{ia}):r\ge0\}\cap S^3$,
$a\in\R/2\pi\Z$. Two distinct pages of $C$ divide $S^3$ into two closed regions,
the \emph{lunes} bounded by these pages. For $C=\{z_2=0\}$ they have the form
$\{(z_1,re^{i\theta}):r\ge0,\ \alpha\le\theta\le\alpha+a\}\cap S^3$ with
$\alpha\in\R$ and $0<a<2\pi$, and $a$ is the \emph{angle} of the lune. A lune of angle $a$ has volume
$\pi a$, its boundary is smooth off $C$, and its perimeter is $4\pi$. A closed
hemisphere is a lune of angle $\pi$. Three distinct pages of $C$ divide $S^3$
into three lunes whose angles add up to $2\pi$, the \emph{three lunes of a
$3$-fan}.

Let $Q=[-1,1]^2$ and
\[
 q(p,t)=\Bigl(\cos\bigl(\tfrac\pi2\|p\|_\infty\bigr)e^{it},\
 \sin\bigl(\tfrac\pi2\|p\|_\infty\bigr)\tfrac p{|p|}\Bigr)\in\C\times\C
 \qquad\bigl((p,t)\in Q\times S^1\bigr),
\]
with $\R^2=\C$ in the second factor and second component $0$ for $p=0$. The map
$q$ is Lipschitz and collapses each circle $\{p\}\times S^1$, $p\in\partial Q$, to
a point. It maps $(\operatorname{int}Q)\times S^1$ homeomorphically onto
$\{z_1\ne0\}$, with nonzero Jacobian almost everywhere. Orient $S^3$ so that $q$
preserves orientation. Identifying $B^2$ with $Q$, a family of loops as in
Section~\ref{sec:widths} gives a continuous map $\Gamma\colon Q\times S^1\to
S^3$ that is constant on $\{p\}\times S^1$ for $p\in\partial Q$. Every such map
is called a \emph{loop family}. Its loops are $\Gamma(p)=\Gamma(p,\cdot)$, and its
degree $\deg\Gamma$ is the degree of the map $\bar\Gamma$ of $S^3$ with
$\Gamma=\bar\Gamma\circ q$.

The next lemma represents the current of a Lipschitz sphere as the boundary of an
integer-valued $BV$ function. In part (d) the density of $Du$ with respect to
$\mathcal H^2$ is the sum of the oriented unit normals of the sheets. So part (d)
bounds every anisotropic total variation $\int f(y,n_u)\,d|Du|$ by the
corresponding anisotropic area of the sphere counted with multiplicity.

\begin{lemma}\label{lem:bounding}
\begin{enumerate}[label=\textup{(\alph*)}]
\item For $u\in L^1(S^3)$, $\mathbf M(\partial(u[\![S^3]\!]))=|Du|(S^3)$, and if
this is finite, then $\partial(u[\![S^3]\!])(\iota_X\mathrm{vol})
=-\int\langle X,dDu\rangle$ for every $C^1$ vector field $X$ on $S^3$.
\item Let $N$ be a compact oriented Riemannian $3$-manifold, possibly with
boundary and corners, and $F\colon N\to S^3$ Lipschitz. Then $F_\#[\![N]\!]=w[\![S^3]\!]$
with $w$ integer-valued and $\|w\|_{L^1}\le\int_N\operatorname{Jac}F$. If
$F_\#[\![\partial N]\!]$ has finite mass, then $w\in BV(S^3;\Z)$.
\item Let $s\colon S^2\to S^3$ be Lipschitz. There is $u\in BV(S^3;\Z)$, unique
up to an additive integer, with $\partial(u[\![S^3]\!])=s_\#[\![S^2]\!]$.
\item In \textup{(c)}, for $z\in S^2$ with $\operatorname{Jac}s(z)>0$ let
$n_s(z)$ be the unit normal of $ds_z(T_zS^2)$ in $T_{s(z)}S^3$ for which
$(n_s(z),ds_ze_1,ds_ze_2)$ is negatively oriented when $(e_1,e_2)$ is a positive
basis of $T_zS^2$. Then
\[
 Du=\Bigl(\sum_{z\in s^{-1}(y)}n_s(z)\Bigr)\,d\mathcal H^2(y),
\]
the sum taken over the points with $\operatorname{Jac}s(z)>0$. Consequently, if
$f\colon TS^3\to[0,\infty)$ is continuous, and convex and positively
homogeneous of degree one on each tangent space, then
\[
 \int f(y,n_u(y))\,d|Du|(y)\le
 \int_{S^2}f(s(z),n_s(z))\operatorname{Jac}s(z)\,dA(z).
\]
In particular $|Du|(S^3)=\mathbf M(s_\#[\![S^2]\!])\le\Area(s)$.
\item For $u\in BV(S^3;\Z)$ and $k\in\Z$ the set $\{u+k\ \text{odd}\}$ has
finite perimeter, and $|D\one_{\{u+k\ \text{odd}\}}|\le|Du|$ as measures.
\end{enumerate}
\end{lemma}

\begin{proof}
(a) For a $2$-form $f$, $\partial(u[\![S^3]\!])(f)=\int_{S^3}u\,df$ depends only
on the restriction of $f$ to $S^3$. Every $2$-form $f$ on $S^3$ extends to $\R^4$
without increasing its comass, as $e(|x|)\pr^*f$ with
$\pr(x)=x/|x|$, $e\in C^\infty_c((0,\infty))$, $0\le e(t)\le t^2$ and $e(1)=1$.
Every $2$-form on $S^3$ is $\iota_X\mathrm{vol}$ for a vector field $X$, with
comass $|X|$ at every point, and $d(\iota_X\mathrm{vol})=\operatorname{div}X\,
\mathrm{vol}$. So the mass is the supremum of $\int u\operatorname{div}X$ over
$|X|\le1$, and the formula is the definition of $Du$.

(b) By the area formula \cite[3.2.3, 4.1.30]{Federer},
$F_\#[\![N]\!](f)=\int_NF^*f=\int_{S^3}w\,f$ for $3$-forms $f$, with
$w(y)=\sum_{z\in F^{-1}(y)}\operatorname{sign}\det dF_z$, the sum over the points
where $dF_z$ exists and is invertible, and $\int|w|\le\int_N\operatorname{Jac}F$.
Since $\partial F_\#[\![N]\!]=F_\#[\![\partial N]\!]$, (a) gives the last claim.

(c) The image of $s$ is $\mathcal H^3$-null, so it misses a point $y_0$. With the
stereographic projection $\mathrm{st}\colon S^3\setminus\{y_0\}\to\R^3$, the map
$k(r,z)=\mathrm{st}^{-1}(r\,\mathrm{st}(s(z)))$ is Lipschitz on $[0,1]\times S^2$,
$k(1,\cdot)=s$ and $k(0,\cdot)$ is constant. So
$\partial k_\#[\![[0,1]\times S^2]\!]=s_\#[\![S^2]\!]$, which has mass at most
$\Area(s)$, and (b) gives $u$. If $\partial(u'[\![S^3]\!])=s_\#[\![S^2]\!]$ as
well, then $|D(u-u')|(S^3)=0$ by (a), and $u-u'$ is constant.

(d) Let $X$ be a $C^1$ vector field and $(e_1,e_2)$ a positive orthonormal frame.
By (a),
\[
 -\int\langle X,dDu\rangle=s_\#[\![S^2]\!](\iota_X\mathrm{vol})
 =\int_{S^2}\mathrm{vol}(X\circ s,ds\,e_1,ds\,e_2)\,dA
 =-\int_{S^2}\langle X\circ s,n_s\rangle\operatorname{Jac}s\,dA .
\]
By the area formula \cite[3.2.3, 3.2.5]{Federer} the last expression is
$-\int\langle X(y),\sum_{z\in s^{-1}(y)}n_s(z)\rangle\,d\mathcal H^2(y)$, and $C^1$
fields are dense in the continuous ones. So $|Du|=|\sum n_s|\,\mathcal H^2$ and
$n_u=\sum n_s/|\sum n_s|$ where $\sum n_s\ne0$. The inequality follows from
$f(y,\sum n_s(z))\le\sum f(y,n_s(z))$ and the area formula for the function
$f(s,n_s)$. For $f(y,v)=|v|$ it gives
$|Du|(S^3)\le\Area(s)$, and $|Du|(S^3)=\mathbf M(s_\#[\![S^2]\!])$ by (a).

(e) Let $k=0$. The function $\one_{\{u\ \text{odd}\}}$ is the $L^1$-limit of
$v_N=\sum_{|j|\le N}(\one_{\{u\ge2j+1\}}-\one_{\{u\ge2j+2\}})$, and
$|Dv_N|(W)\le\sum_k|D\one_{\{u\ge k\}}|(W)=|Du|(W)$ for open $W$ by the coarea
formula. Lower semicontinuity gives $|D\one_{\{u\ \text{odd}\}}|(W)\le|Du|(W)$
for open, hence for Borel, sets $W$.
\end{proof}

The function $u$ of (c) is a \emph{bounding function} of $s$. The complementary
sets $\{u\ \text{odd}\}$ and $\{u\ \text{even}\}$ are the \emph{odd sets} of $s$.
By (d) and (e) their perimeters are at most $|Du|(S^3)\le\Area(s)$.

The next lemma computes the degree of a loop family from its push-forward.

\begin{lemma}\label{lem:degree}
Let $\Gamma\colon Q\times S^1\to S^3$ be a Lipschitz loop family. Then
$\Gamma_\#[\![Q\times S^1]\!]=(\deg\Gamma)[\![S^3]\!]$.
\end{lemma}

\begin{proof}
Let $G\colon S^3\to S^3$ be smooth with $|G-\bar\Gamma|<\frac12$. By
Lemma~\ref{lem:bounding}(b), applied to $G$ and to $q$, and by Sard's theorem,
$G_\#[\![S^3]\!]=(\deg G)[\![S^3]\!]$ and $q_\#[\![Q\times S^1]\!]=[\![S^3]\!]$;
so $(G\circ q)_\#[\![Q\times S^1]\!]=(\deg\Gamma)[\![S^3]\!]$. With
$\pr(x)=x/|x|$, the map $H(r,\cdot)=\pr((1-r)\Gamma+r\,G\circ q)$ is Lipschitz on
$[0,1]\times Q\times S^1$, and the homotopy formula \cite[4.1.9, 4.1.14]{Federer} gives,
for $T=[\![Q\times S^1]\!]$,
\[
 (G\circ q)_\#T-\Gamma_\#T=\partial H_\#([\![0,1]\!]\times T)
 +H_\#([\![0,1]\!]\times\partial T).
\]
The $4$-current $H_\#([\![0,1]\!]\times T)$ is zero, since $H$ takes values in
$S^3$. On $[0,1]\times\partial Q\times S^1$ the map $H$ does not depend on the
circle variable, so its image is $\mathcal H^3$-null and the last current is
zero.
\end{proof}

We glue the fillings of two loops and an annulus between the loops into a sphere.
The following theorem selects three nearby loops and one odd set for each of the
three switching spheres. These sets cover each point an odd number of times
outside a set of volume at most $\eps$.

\begin{theorem}[switching triangles]\label{thm:switching}
Let $\Gamma\colon Q\times S^1\to S^3$ be a loop family of degree $\pm1$ whose
loops $\Gamma(p)$ are $L$-Lipschitz, and let $\eps>0$. There is $h_0>0$ with the
following property. Let $\mathcal T$ be a triangulation of $Q$ by triangles of
diameter at most $h_0$, and for each vertex $x$ of $\mathcal T$ let
$f_x\colon\bar D\to S^3$ be a Lipschitz disk with $f_x|_{\partial D}=\Gamma(x)$,
constant if $x\in\partial Q$. For an edge $[x,y]$ let
\[
 a_{xy}(r,t)=\pr\bigl((1-r)\Gamma(x,t)+r\,\Gamma(y,t)\bigr)\qquad
 \bigl((r,t)\in[0,1]\times S^1\bigr),
\]
and let $s_{xy}\colon S^2\to S^3$ be the \emph{switching sphere} glued from
$f_x$, $a_{xy}$ and $f_y$ along their boundary loops, oriented so that
$s_{xy\#}[\![S^2]\!]=f_{x\#}[\![D]\!]+a_{xy\#}[\![[0,1]\times S^1]\!]
-f_{y\#}[\![D]\!]$. Then there are a triangle of $\mathcal T$ with vertices
$x_1,x_2,x_3$ and bounding functions $u_{12},u_{23},u_{31}$ of
$s_{x_1x_2},s_{x_2x_3},s_{x_3x_1}$ such that the odd sets
$E_{ij}=\{u_{ij}\ \text{odd}\}$ satisfy:
\begin{enumerate}[label=\textup{(\alph*)}]
\item the loops satisfy $|\Gamma(x_i,t)-\Gamma(x_j,t)|\le\eps$ for all $t$,
and the annuli have $\Area(a_{x_ix_j})\le\eps$;
\item each $E_{ij}$ has volume at most $\pi^2$;
\item the sum $\one_{E_{12}}+\one_{E_{23}}+\one_{E_{31}}$ is odd outside a set
of volume at most $\eps$.
\end{enumerate}
\end{theorem}

\begin{proof}
Let $\delta(h)=\sup\{|\Gamma(p,t)-\Gamma(p',t)|:|p-p'|\le h,\ t\in S^1\}$, which tends
to $0$ with $h$, and choose $h_0$ with $\delta(h_0)\le\min\{\eps,\frac12\}$,
$8\pi L\delta(h_0)\le\eps$ and $8\pi L\delta(h_0)^2\le\eps$. Let $\mathcal T$ and the $f_x$
be given, and write $\delta=\delta(h_0)$. For $p$ in a triangle with vertices
$v_0,v_1,v_2$ and barycentric coordinates $b_i(p)$ put
$\Gamma'(p,t)=\pr(\sum_ib_i(p)\Gamma(v_i,t))$. The convex combination has length
at least $1-\delta\ge\frac12$, so $\Gamma'$ is well defined, Lipschitz, and a loop
family. It equals $a_{xy}$ on each edge. The maps $\pr((1-r)\Gamma(p,t)+
r\sum_ib_i(p)\Gamma(v_i,t))$, $r\in[0,1]$, are well defined by the same estimate
and join $\Gamma$ to $\Gamma'$ through loop families. So $\deg\Gamma'=\pm1$. Since
$\|d\pr\|\le2$ where $|x|\ge\frac12$, the annulus $a_{xy}$ has $|\partial_ra_{xy}|
\le2\delta$ and $|\partial_ta_{xy}|\le2L$, so $\Area(a_{xy})\le8\pi L\delta\le\eps$. On a
triangle $K$, in the barycentric coordinates $(b_1,b_2)$ of the standard
simplex, the derivatives of $\Gamma'$ in $b_1,b_2$ have length at most $2\delta$. So,
by Lemma~\ref{lem:bounding}(b), $T_K=\Gamma'_\#[\![K\times
S^1]\!]=w_K[\![S^3]\!]$ with $w_K$ integer-valued and
$\|w_K\|_{L^1}\le\frac12\cdot2\pi\cdot8L\delta^2\le\eps$.

For an oriented edge $e=[x,y]$ let $Z_e=s_{xy\#}[\![S^2]\!]$. If $e\subset\partial
Q$, then $Z_e=0$, since $f_x$, $f_y$ and $a_{xy}$ have images of dimension at
most one. Choose bounding functions $u_e$ of $s_{xy}$ with $u_e=0$ if $e\subset\partial Q$
and with $|\{u_e\ \text{odd}\}|\le\pi^2$. This is possible, as the odd sets of $u_e$
and $u_e+1$ are complementary. For the reversed edge, whose cycle is $-Z_e$, take
$u_{-e}=-u_e$. For a positively oriented triangle $K=[x,y,z]$ the
disk terms cancel in $Z_{xy}+Z_{yz}+Z_{zx}$, which is therefore
$\Gamma'_\#[\![\partial K\times S^1]\!]=\partial T_K$. By
Lemma~\ref{lem:bounding}(a), the function $u_{xy}+u_{yz}+u_{zx}-w_K$ has
zero variation, so
\[
 u_{xy}+u_{yz}+u_{zx}=w_K+j_K\qquad\text{almost everywhere}
\]
for an integer $j_K$. Summing over all
positively oriented triangles, the terms of interior edges cancel, boundary
edges contribute zero, and $\sum_K w_K=\pm1$ by
Lemma~\ref{lem:degree}, as $\sum_K T_K=\Gamma'_\#[\![Q\times
S^1]\!]$. So $\sum_K j_K=\mp1$, and $j_K$ is odd for some triangle
$K=[x_1,x_2,x_3]$. On $\{w_K=0\}$, reduction modulo $2$ gives (c), and
$|\{w_K\ne0\}|\le\|w_K\|_{L^1}\le\eps$. The bounds (a) and (b) hold by
construction.
\end{proof}

\begin{remark}
The proof uses only the degree of $\Gamma$, its uniform continuity, the Lipschitz
bound of the loops and the given disks. No bound on the areas of the disks
enters. If the disks have area at most $A+\eps$, the switching spheres have area
at most $2A+3\eps$. By Lemma~\ref{lem:bounding} the three odd sets then have
perimeter at most $2A+3\eps$, and they cover $S^3$ an odd number of times up to
volume $\eps$.
\end{remark}

The next proposition passes to the limit along a sequence of switching triangles.
If the variations are uniformly bounded, the bounding functions converge along a
subsequence after adding even integers.

\begin{proposition}[limits]\label{prop:limits}
For $i=1,2,3$ and $n\ge1$ let $u_i^n\in BV(S^3;\Z)$ with
$\sup_n|Du_i^n|(S^3)<\infty$, such that $E_i^n=\{u_i^n\ \text{odd}\}$ has
$|E_i^n|\le\pi^2$ and $\one_{E^n_1}+\one_{E^n_2}+\one_{E^n_3}$ is odd outside a
set of volume tending to $0$. After passing to a subsequence and adding even
integers to the $u^n_i$, there are $u_i\in BV(S^3;\Z)$ such that $u_i^n\to u_i$
in $L^1$ and almost everywhere, $Du^n_i\to Du_i$ weakly-*, and
$|Du_i|(S^3)\le\liminf_n|Du_i^n|(S^3)$. The sets $E_i=\{u_i\ \text{odd}\}$
satisfy $\one_{E_i^n}\to\one_{E_i}$ in $L^1$, $|E_i|\le\pi^2$,
$P(E_i)\le|Du_i|(S^3)$, and $\one_{E_1}+\one_{E_2}+\one_{E_3}$ is odd almost
everywhere. So almost every point lies in one or in all three of the $E_i$,
and $|E_1|+|E_2|+|E_3|\ge2\pi^2$.
\end{proposition}

\begin{proof}
Let $u=u^n_i$, and let $l$ be the least integer with $|\{u>l\}|\le\pi^2$. Adding
an even integer changes neither $E^n_i$ nor $Du$, so let $l\in\{0,1\}$.
Then $|\{u\ge2\}|\le\pi^2$ and $|\{u\le-1\}|<\pi^2$. By the layer-cake formula
\[
 \|u\|_{L^1}=\sum_{j\ge1}\bigl(|\{u\ge j\}|+|\{u\le-j\}|\bigr),
\]
and every set in this sum except $\{u\ge1\}$ has volume at most $\pi^2$, hence
at most $c\,P^{3/2}$ by \eqref{eq:schmidt}. Since $P(\{u\le-j\})=P(\{u\ge1-j\})$
and $\sum P^{3/2}\le(\sum P)^{3/2}$, the coarea formula gives
$\|u\|_{L^1}\le2\pi^2+c\,|Du|(S^3)^{3/2}$. So $(u^n_i)_n$ is bounded in
$BV(S^3)$, and a subsequence converges in $L^1$ and almost everywhere to an
integer-valued $u_i$, with the stated convergence of $Du^n_i$ and
semicontinuity. Integer-valued sequences converging almost everywhere are
eventually constant at almost every point, so $\one_{E_i^n}\to\one_{E_i}$
almost everywhere and in $L^1$, $|E_i|\le\pi^2$, and the parity condition holds
almost everywhere in the limit. Lemma~\ref{lem:bounding}(e) gives
$P(E_i)\le|Du_i|(S^3)$, and $\sum_i\one_{E_i}\ge1$ almost everywhere gives the
volume bound.
\end{proof}

If $|Du|(S^3)\le4\pi$ and the odd set is a lune, then the function is constant on
each side of its boundary, and the two constants differ by one. This is the next
proposition.

\begin{proposition}[lune criterion]\label{prop:lune}
Let $u\in BV(S^3;\Z)$ with $|Du|(S^3)\le4\pi$, $k\in\Z$, and suppose that
$E=\{u+k\ \text{odd}\}$ is almost everywhere equal to a lune $L$, for instance a
closed hemisphere. Then $P(E)=4\pi$, $|Du|=|D\one_E|=\mathcal
H^2\llcorner\partial L$, and $u=b\pm\one_L$ almost everywhere for an integer
$b$.
\end{proposition}

\begin{proof}
Let $L$ be the lune between the pages $P_1,P_2$ of $C$. Its boundary is smooth
off $C$ and $\mathcal H^2(C)=0$, so the divergence theorem gives
$|D\one_E|=|D\one_L|=\mathcal H^2\llcorner(P_1\cup P_2)$ and $P(E)=4\pi$. By
Lemma~\ref{lem:bounding}(e), $|D\one_E|\le|Du|$, and $|Du|(S^3)\le4\pi$; hence
$|Du|=|D\one_E|$. So $u$ is almost everywhere constant on each of the two
components of $S^3\setminus(P_1\cup P_2)$, equal to $b_1$ on $L$ and to $b$ on
its complement. Then $|Du|=|b_1-b|\,|D\one_L|$, so $b_1-b=\pm1$.
\end{proof}

\begin{remark}\label{rem:limits}
Consider switching triangles of Theorem~\ref{thm:switching} with $\eps_n\to0$ and
switching spheres of area at most $4\pi+o(1)$. Proposition~\ref{prop:limits}
gives limits of the odd sets. They have perimeter at most $4\pi$ and volume at
most $\pi^2$, and they cover almost every point an odd number of times. A limit of volume
$\pi^2$ is a closed hemisphere by Schmidt's inequality. If all three have volume
less than $\pi^2$, they are the lunes of a $3$-fan by Theorem~D. If a limit odd
set is a lune or a closed hemisphere, then by Proposition~\ref{prop:lune} its
bounding function jumps by one across the boundary and nowhere else.
\end{remark}

\begin{remark}
In Theorem~D the volumes add up to at least $2\pi^2$, so each is at least
$2(\pi^2-\max_j|E_j|)>0$. The bound $|E_i|<\pi^2$ cannot be weakened to
$|E_i|\le\pi^2$: two complementary closed hemispheres and the empty set satisfy
the other hypotheses.
\end{remark}

\section{Three-set rigidity}\label{sec:threeset}

In this section we prove Theorem~D. Lemma~\ref{lem:threecell} assigns the points
covered three times to one of the sets, which reduces the theorem to partitions of
$S^3$ into three sets. Proposition~\ref{prop:cone} relates the perimeters of such a
partition to the Gaussian perimeters of the cones over the three sets. The main
tool is a Gaussian double bubble theorem with unequal interface weights
(Theorem~\ref{thm:wgdb}). Its inequality and its equality case, applied to the
cones, give the sharp weighted perimeter bound of Corollary~\ref{cor:fan}, from
which Theorem~D follows.

We first fix notation. In $\R^k$, $|X|$ is the Lebesgue measure of $X$ and $B_R$
the open ball of radius $R$ about $0$. For subsets of reduced boundaries and of
hyperplanes, null sets are $\mathcal H^{k-1}$-null sets. For a set $X$ of locally
finite perimeter in $\R^k$ or in $S^3$ (then $k=3$), $\partial^*X$ is its reduced
boundary, and $|D\one_X|=\mathcal H^{k-1}\llcorner\partial^*X$
\cite[Theorem~15.9]{Maggi}. The set $\partial^*X\cap W$, for $W$ open, depends
only on $X\cap W$. If it is $\mathcal H^{k-1}$-null and $W$ is connected, then
$X\cap W$ is null or of full measure in $W$. Let $X_1,\dots,X_l$ be pairwise disjoint sets of
locally finite perimeter covering $S^3$ or $\R^k$ up to a null set. Then
$\mathcal H^{k-1}$-almost every point of $\partial^*X_i$ lies in $\partial^*X_j$
for exactly one $j\ne i$, and there the measure-theoretic normals of $X_i$ and
$X_j$ are opposite \cite[Theorem~4.17]{AFP}, \cite[Theorems~15.5
and~16.3]{Maggi}. On $S^3$ this gives, with
$p_{ij}=\mathcal H^2(\partial^*X_i\cap\partial^*X_j)$,
\begin{equation}\label{eq:partition}
 P(X_i)=\sum_{j\ne i}p_{ij},\qquad P(X_i\cup X_j)=P(X_i)+P(X_j)-2p_{ij}\quad(i\ne j).
\end{equation}

The first lemma assigns the points covered three times by the sets of Theorem~D to
one of the cells, so that it suffices to treat partitions of $S^3$.

\begin{lemma}\label{lem:threecell}
Let $E_1,E_2,E_3$ be as in Theorem~D. After renumbering the
$E_i$, there is a partition $(Y_1,Y_2,Y_3)$ of $S^3$, up to null sets, into sets
of finite perimeter with $|Y_i|<\pi^2$ and $P(Y_i)\le4\pi$, such that $Y_i=E_i$
up to null sets for every $i$ if $P(Y_1)=P(Y_2)=P(Y_3)=4\pi$.
\end{lemma}

\begin{proof}
Let $X_i$ be the set of points of $E_i$ in no other $E_j$, and
$X_4=E_1\cap E_2\cap E_3$. By the parity condition $X_1,\dots,X_4$ partition
$S^3$ up to null sets, and $E_i=X_i\cup X_4$ up to null sets. Renumber the $E_i$ so that
$p_{41}\ge p_{42},p_{43}$, and put $Y_1=E_1$, $Y_2=X_2$, $Y_3=X_3$. By
\eqref{eq:partition}, $P(E_2)=P(X_2)+P(X_4)-2p_{24}=p_{21}+p_{23}+p_{41}+p_{43}$,
so
\[
 P(Y_2)=P(E_2)+p_{24}-p_{41}-p_{43}\le P(E_2)-p_{43},\qquad\text{and likewise}\qquad
 P(Y_3)\le P(E_3)-p_{42}.
\]
No volume increases. If all $P(Y_i)=4\pi$, then $p_{43}=p_{42}=0$, and then the
identity for $P(Y_2)$ gives $p_{41}\le p_{24}=0$. So $P(X_4)=0$, and $X_4$ is
null, as $|X_4|<2\pi^2$.
\end{proof}

We next study clusters in $\R^n$ with the Gaussian measure. Let $n\ge2$. The
standard Gaussian measure on $\R^n$
and its density are both denoted $\gamma$, $\gamma(x)=(2\pi)^{-n/2}e^{-|x|^2/2}$.
The index in $\gamma_1,\gamma_2$ indicates the dimension when it differs from
$n$. Let $N(s)=\gamma_1((-\infty,s))$, $I_\gamma(v)=\gamma_1(N^{-1}(v))$ for
$0<v<1$, $I_\gamma(0)=I_\gamma(1)=0$, and $\gamma^{n-1}=\gamma\,\mathcal
H^{n-1}$. For a set $E\subset\R^n$ of locally finite perimeter, $n_E$ is the
outer unit normal on $\partial^*E$, so that $\int_E\operatorname{div}Y=
\int_{\partial^*E}\langle Y,n_E\rangle\,d\mathcal H^{n-1}$ for $Y\in
C^1_c(\R^n;\R^n)$. The Gaussian perimeter
$P_\gamma(E)=\sup\{\int_E\operatorname{div}(\gamma Y):Y\in C^1_c(\R^n;\R^n),\
|Y|\le1\}$ is finite only if $E$ has locally finite perimeter, as $\gamma$ is
bounded below on bounded sets. Then $P_\gamma(E)=\gamma^{n-1}(\partial^*E)$.
By the Gaussian isoperimetric inequality \cite{Borell}, \cite{ST},
$P_\gamma(E)\ge I_\gamma(\gamma(E))$, with equality for half-spaces. In this
form it follows from Bobkov's inequality $I_\gamma(\int f\,d\gamma)-\int
I_\gamma(f)\,d\gamma\le\int|\nabla f|\,d\gamma$ for locally Lipschitz
$f\colon\R^n\to[0,1]$ \cite[Corollary~2]{Bobkov}. If $E\subset B_R$, let
$f=\one_E*k_\delta$ be the mollifications, with $k_\delta\ge0$ smooth, of integral
one and supported in $B_\delta$, $\delta<1$. They satisfy $\int|\nabla
f|\,d\gamma\le e^{\delta(R+1)}P_\gamma(E)$, as $\gamma(x)\le
e^{\delta(R+1)}\gamma(y)$ for $|x-y|<\delta$, $|y|\le R$. They converge to $\one_E$
almost everywhere as $\delta\to0$. In general, for almost every $R$,
$P_\gamma(E\cap B_R)\le P_\gamma(E)+\gamma^{n-1}(E\cap\partial B_R)$
\cite[Lemma~15.12]{Maggi}. Since $\int_0^\infty\gamma^{n-1}(E\cap\partial
B_R)\,dR=\gamma(E)$, the bounded case applies along radii $R\to\infty$ with
$\gamma^{n-1}(E\cap\partial B_R)\to0$.

A \emph{cluster} is a triple $U=(U_1,U_2,U_3)$ of sets of locally finite
perimeter that are pairwise disjoint up to null sets and cover $\R^n$ up to a null
set. Clusters that agree up to null sets are identified. The vector
$\gamma(U)=(\gamma(U_i))_i$ lies in the simplex
$\mathcal V=\{v\in\R^3:v_i\ge0,\ \sum v_i=1\}$. Its interior is
$\mathcal V^\circ$, and its tangent space is $\R^3_0=\{y\in\R^3:\sum y_i=0\}$.
Let $S_{ij}=\partial^*U_i\cap\partial^*U_j$, $a_{ij}=\gamma^{n-1}(S_{ij})$,
and let $n_{ij}=n_{U_i}=-n_{U_j}$ on $S_{ij}$. Throughout
$\{i,j,k\}=\{1,2,3\}$. Fix $w_1,w_2,w_3>0$ and put $w_{ij}=w_i+w_j$. Following
Lawlor \cite{Lawlor}, we call $w_{ij}$ the \emph{weight} of the interface
$S_{ij}$. The weights satisfy the strict triangle inequality, and all weights
with the strict triangle inequality arise in this way. The weighted perimeter and the weighted profile are
\[
 F_w(U)=\sum_iw_iP_\gamma(U_i)=\sum_{i<j}w_{ij}a_{ij},\qquad
 I_w(v)=\inf\{F_w(U):\gamma(U)=v\}\quad(v\in\mathcal V),
\]
and $U$ is a \emph{minimizer} if $F_w(U)=I_w(\gamma(U))$. Quadratic forms on
$\R^3_0$ are identified with symmetric operators of $\R^3_0$, ordered as forms.
For $b=(b_{ij})_{i<j}$ with $b_{ij}\ge0$ the operator
$L_b=\sum_{i<j}b_{ij}(e_i-e_j)(e_i-e_j)^T$ of $\R^3$ vanishes on $(1,1,1)$ and
maps $\R^3_0$ to itself. We write $L_b^{-1}$ for the inverse of its restriction to
$\R^3_0$, when it exists, and $(a/w)_{ij}=a_{ij}/w_{ij}$.

Choose $c_1,c_2,c_3\in\R^2$ with $|c_i-c_j|=w_{ij}$, let
$A$ be the $3\times2$ matrix with rows $c_i^T$, and let
$T_i=\{z\in\R^2:\langle c_i-c_j,z\rangle>0\ \text{for } j\ne i\}$. These are three
open sectors with vertex $0$. The angle of $T_k$ is $\pi$ minus the angle of the
triangle $c_1c_2c_3$ at $c_k$. The unit vector $n_{ij}=(c_j-c_i)/w_{ij}$ is
normal to the ray $\bar T_i\cap\bar T_j$ and points from $T_i$ to $T_j$. Let
$d_{ij}=d_{ji}$ be the unit direction of this ray. For $(i,j)$ in cyclic order
$d_{ij}$ is $n_{ij}$ turned by a fixed right angle, so $\sum w_{ij}n_{ij}=0$
(cyclic sum) gives \emph{Young's law} $\sum_{i<j}w_{ij}d_{ij}=0$. A
\emph{weighted tripod cluster} is a cluster
$U_i=\{x:\langle c'_i,x\rangle-l_i>\langle c'_j,x\rangle-l_j\ \text{for } j\ne i\}$
with $c'_i\in\R^n$, $|c'_i-c'_j|=w_{ij}$ and $l\in\R^3$. As the triangles
$c'_1c'_2c'_3$ and $c_1c_2c_3$ are congruent, there is an orthogonal map $O$ of
$\R^n$ with $O(c'_i-c'_j)=(c_i-c_j,0)\in\R^2\times\R^{n-2}$. There is also
$x\in\R^2$ with $\langle c_i-c_j,x\rangle=l_i-l_j$ for all $i,j$. Then
$OU_i=(x+T_i)\times\R^{n-2}$. For
$x\in\R^2$ let
\[
 V(x)=\bigl(\gamma_2(x+T_i)\bigr)_i,\qquad
 \tilde a_{ij}(x)=\int_{x+(\bar T_i\cap\bar T_j)}\gamma_2\,d\mathcal H^1,\qquad
 \tilde F(x)=\sum_{i<j}w_{ij}\tilde a_{ij}(x).
\]

The next lemma records the Gauss--Green formula for half-spaces and wedges. It also
gives the volumes and the interfaces of the clusters formed by the sectors.

\begin{lemma}\label{lem:wedge}
Let $a,b\in\R^n$ be linearly independent unit vectors and $a_0,b_0\in\R$.
\begin{enumerate}[label=\textup{(\alph*)}]
\item The half-space $\{\langle a,x\rangle<a_0\}$ has reduced boundary
$\{\langle a,x\rangle=a_0\}$ and outer normal $a$.
\item The wedge $E=\{\langle a,x\rangle<a_0,\ \langle b,x\rangle<b_0\}$
satisfies $\int_E\operatorname{div}Y=\int_{F_a}\langle Y,a\rangle\,d\mathcal
H^{n-1}+\int_{F_b}\langle Y,b\rangle\,d\mathcal H^{n-1}$ for $Y\in
C^1_c(\R^n;\R^n)$, where $F_a=\{\langle a,x\rangle=a_0,\ \langle
b,x\rangle<b_0\}$ and $F_b=\{\langle b,x\rangle=b_0,\ \langle a,x\rangle<a_0\}$.
So $\partial^*E=F_a\cup F_b$ up to a null set, with outer normals $a,b$.
\item For $x\in\R^2$ the sets $(x+T_i)\times\R^{n-2}$ form a cluster $U$ with
$S_{ij}=(x+((\bar T_i\cap\bar T_j)\setminus\{0\}))\times\R^{n-2}$ up to an
$\mathcal H^{n-1}$-null set,
$\gamma(U)=V(x)$, $a_{ij}=\tilde a_{ij}(x)$ and $F_w(U)=\tilde F(x)$.
\end{enumerate}
\end{lemma}

\begin{proof}
(b) Write $Y=\langle Y,a\rangle a^*+\langle Y,b\rangle b^*+Y'$ with
$a^*,b^*\in\operatorname{span}\{a,b\}$ dual to $a,b$ and $Y'\perp a,b$. As $E$ is
invariant under translations orthogonal to $a,b$, Fubini's theorem gives
$\int_E\operatorname{div}Y'=0$. Integrating $\partial_{a^*}\langle Y,a\rangle$
along the lines $y+ra^*$, $y\in\{\langle a,y\rangle=a_0\}$, on which
$\langle a,\cdot\rangle$ increases with unit speed and
$\langle b,\cdot\rangle$ is constant, gives the term of $F_a$; likewise for $b$.
Part (a) follows in the same way. (c) The cell $x+T_i$ is the wedge
$\{\langle n_{ij},z-x\rangle<0,\ \langle n_{ik},z-x\rangle<0\}$, whose face
toward $x+T_j$ is the open ray $x+((\bar T_i\cap\bar T_j)\setminus\{0\})$. Apply (b)
and use $\gamma(z',z'')=\gamma_2(z')\gamma_{n-2}(z'')$.
\end{proof}

The next lemma says that clusters of uniformly bounded weighted perimeter are
sequentially compact. It also gives the existence of minimizers and the boundary
values of the profile.

\begin{lemma}\label{lem:exist}
\begin{enumerate}[label=\textup{(\alph*)}]
\item If $U^h$ are clusters with $\sup_hF_w(U^h)<\infty$, a subsequence converges
cell by cell in $L^1_{\mathrm{loc}}(\R^n)$ and in $L^1(\gamma)$ to a cluster $U$
with $\gamma(U)=\lim_h\gamma(U^h)$ and $F_w(U)\le\liminf_hF_w(U^h)$.
\item The profile $I_w$ is bounded and lower semicontinuous on $\mathcal V$,
and every $v\in\mathcal V$ is $\gamma(U)$ for a minimizer $U$.
\item If $\gamma(U)\in\mathcal V^\circ$ and $F_w(U)<\infty$, then at least two
of the $a_{ij}$ are positive and $L_{a/w}$ is positive definite on $\R^3_0$.
\item If $v\in\partial\mathcal V$, $v_k=0<v_i,v_j$, then
$I_w(v)=w_{ij}I_\gamma(v_i)$; $I_w=0$ at the vertices; and
$I_w(v)\ge\min_lw_l\sum_iI_\gamma(v_i)$ for all $v\in\mathcal V$.
\end{enumerate}
\end{lemma}

\begin{proof}
(a) Each $P_\gamma(U^h_i)$ is at most $\sup_hF_w(U^h)/\min_lw_l$, so the
perimeters of $U^h_i$ in every ball are bounded. By compactness in
$BV_{\mathrm{loc}}$ \cite[Theorem~3.23]{AFP} and a diagonal argument, a
subsequence has $\one_{U^h_i}\to f_i$ in $L^1_{\mathrm{loc}}$ and almost
everywhere. Then $f_i\in\{0,1\}$ and $\sum f_i=1$ almost everywhere, and the
sets $\{f_i=1\}$, corrected on null sets, form a partition $U$. Dominated
convergence gives convergence in $L^1(\gamma)$ and of the volumes. For fixed
$Y$, $E\mapsto\int_E\operatorname{div}(\gamma Y)$ is continuous in
$L^1_{\mathrm{loc}}$, so $P_\gamma$ is lower semicontinuous. Hence
$P_\gamma(U_i)<\infty$, $U$ is a cluster, and $F_w(U)\le\liminf F_w(U^h)$ as
$w_i>0$.

(b) For $v\in\mathcal V$ the slabs $\{x_1<s_1\}$, $\{s_1<x_1<s_2\}$, $\{x_1>s_2\}$
with $-\infty\le s_1\le s_2\le\infty$, $N(s_1)=v_1$, $N(s_2)=v_1+v_2$ have volumes
$v$; put $\gamma_1(\pm\infty)=0$. Their interfaces lie in
the hyperplanes $\{x_1=s_1\}$ and $\{x_1=s_2\}$. The perimeter of the middle slab is
at most the sum of the perimeters of the two half-spaces of which it is the
difference. So Lemma~\ref{lem:wedge}(a) gives $F_w\le
w_{12}\gamma_1(s_1)+w_{23}\gamma_1(s_2)\le\frac2{\sqrt{2\pi}}\max w_{ij}$. Let $v^h\to v$, along a subsequence realizing
$\liminf I_w(v^h)$, and $\gamma(U^h)=v^h$ with $F_w(U^h)\le I_w(v^h)+\frac1h$.
Then (a) gives a cluster $U$ with $\gamma(U)=v$ and
$I_w(v)\le F_w(U)\le\liminf I_w(v^h)$; for $v^h=v$, $U$ is a minimizer.

(c) If two of the $a_{ij}$ vanish, they share an index $k$, and
$P_\gamma(U_k)=a_{ik}+a_{jk}=0$, whereas $P_\gamma(U_k)\ge I_\gamma(v_k)>0$. Any
two of the $e_i-e_j$ span $\R^3_0$, and
$\langle L_{a/w}y,y\rangle=\sum_{i<j}(a_{ij}/w_{ij})(y_i-y_j)^2$.

(d) If $v_k=0$, then $U_k$ is null, $U_j=\R^n\setminus U_i$ up to a null set,
and $F_w(U)=w_{ij}P_\gamma(U_i)\ge w_{ij}I_\gamma(v_i)$, with equality for
half-spaces by Lemma~\ref{lem:wedge}(a). At a vertex the cluster
$(\R^n,\emptyset,\emptyset)$ has $F_w=0$. The last bound is the isoperimetric
inequality for each cell.
\end{proof}

We vary clusters by the maps $R_t$ of the next lemma. It gives the first two
derivatives of the volumes and of the perimeters.

\begin{lemma}\label{lem:rankone}
Let $e\in\R^n$, $f\in C^\infty_c(\R^n)$ with support $K$,
$|t|\,|e|\sup|\nabla f|\le\frac12$, and $R_t(x)=x+tf(x)e$.
\begin{enumerate}[label=\textup{(\alph*)}]
\item The map $R_t$ is a diffeomorphism of $\R^n$, the identity off $K$, with
$\det DR_t=1+t\partial_ef\ge\frac12$. If $E$ has locally finite perimeter, so
has $R_tE$, and with $W_t=(1+t\partial_ef)n_E-t\langle e,n_E\rangle\nabla f$,
\[
 \gamma(R_tE)=\int_E(\gamma\circ R_t)(1+t\partial_ef)\,dx,\qquad
 P_\gamma(R_tE)=\int_{\partial^*E}(\gamma\circ R_t)\,|W_t|\,d\mathcal H^{n-1}.
\]
If $U$ is a cluster, so is $R_tU=(R_tU_i)_i$.
\item Let $P_\gamma(E)<\infty$, $m=\gamma^{n-1}\llcorner\partial^*E$,
$p=\langle e,n_E\rangle$, and for $h\in C^1(\R^n)$ let
$\nabla^\top h=\nabla h-\langle\nabla h,n_E\rangle n_E$ and
$\ell(h)=\langle\nabla^\top h,e\rangle$. For $h\in C^\infty_c(\R^n)$ put
\begin{gather*}
 \delta V_E(h)=\int hp\,dm,\qquad
 \delta P_E(h)=\int\bigl(\ell(h)-h\langle e,x\rangle\bigr)\,dm,\qquad
 \delta^2V_E(f)=-\delta V_E\bigl(f^2\langle e,x\rangle\bigr),\\
 \delta^2P_E(f)=\int\Bigl(p^2|\nabla^\top f|^2-2f\langle e,x\rangle\ell(f)
 +f^2\bigl(\langle e,x\rangle^2-|e|^2\bigr)\Bigr)\,dm .
\end{gather*}
There are $c,t_0>0$, depending only on $f$ and $e$, such that for $|t|\le t_0$
\begin{align*}
 \bigl|\gamma(R_tE)-\gamma(E)-t\,\delta V_E(f)-\tfrac12t^2\delta^2V_E(f)\bigr|
 &\le c|t|^3,\\
 \bigl|P_\gamma(R_tE)-P_\gamma(E)-t\,\delta P_E(f)-\tfrac12t^2\delta^2P_E(f)\bigr|
 &\le c|t|^3m(K).
\end{align*}
\end{enumerate}
\end{lemma}

\begin{proof}
(a) If $R_t(x)=R_t(y)$, then $|x-y|=|t||e||f(x)-f(y)|\le\frac12|x-y|$, so
$R_t$ is injective. It maps each line $x+\R e$ onto itself, as $r\mapsto
r+tf(x+re)$ is continuous and equal to $r$ for large $|r|$, and
$DR_t=I+te\,\nabla f^T$ has determinant $1+t\partial_ef$. For $Y\in
C^1_c(\R^n;\R^n)$ the field $S=\operatorname{adj}(DR_t)\,(Y\circ
R_t)=(1+t\partial_ef)\,Y\circ R_t-t\langle\nabla f,Y\circ R_t\rangle e$ satisfies
$\operatorname{div}S=(1+t\partial_ef)(\operatorname{div}Y)\circ R_t$ (Piola's
identity; directly, the terms containing $\nabla\partial_ef$ and $\partial_eY$
cancel), and $\langle S,n_E\rangle=\langle Y\circ R_t,W_t\rangle$. By the change of
variables and the Gauss--Green formulas,
\[
 \int_{R_tE}\operatorname{div}Y=\int_E\operatorname{div}S
 =\int_{\partial^*E}\langle Y\circ R_t,W_t\rangle\,d\mathcal H^{n-1}.
\]
So $-D\one_{R_tE}$ is the push-forward of $W_t\,\mathcal H^{n-1}\llcorner\partial^*E$
under the homeomorphism $R_t$, and $|D\one_{R_tE}|$ that of
$|W_t|\,\mathcal H^{n-1}\llcorner\partial^*E$; integrating $\gamma$ gives
$P_\gamma(R_tE)$. The volume formula is the change of variables formula, and
$R_t$ preserves disjointness and null sets.

(b) Since $\partial_ef=\langle\nabla f,n_E\rangle p+\ell(f)$,
$W_t=(1+t\ell(f))n_E-tp\nabla^\top f$ and $|W_t|^2=(1+t\ell(f))^2+t^2p^2|\nabla^\top
f|^2$; and $\gamma\circ R_t=\gamma\exp(-tf\langle e,x\rangle-\frac12t^2f^2|e|^2)$.
On $K$ all these functions are bounded in terms of $f$ and $e$; off $K$,
$R_t=\mathrm{id}$ and $|W_t|=1$. Taylor's formula gives, uniformly on $K$,
\begin{align*}
 \tfrac{\gamma\circ R_t}\gamma(1+t\partial_ef)&=1+t\bigl(\partial_ef-f\langle
 e,x\rangle\bigr)+\tfrac12t^2\bigl(f^2\langle e,x\rangle^2-f^2|e|^2
 -2f\langle e,x\rangle\partial_ef\bigr)+O(t^3),\\
 \tfrac{\gamma\circ R_t}\gamma|W_t|&=1+t\bigl(\ell(f)-f\langle e,x\rangle\bigr)
 +\tfrac12t^2\bigl(p^2|\nabla^\top f|^2-2f\langle e,x\rangle\ell(f)\\
 &\qquad+f^2(\langle e,x\rangle^2-|e|^2)\bigr)+O(t^3).
\end{align*}
Integrating the second line against $m$ gives the expansion of $P_\gamma$. In the
first line, the coefficients of $t$ and of $\frac12t^2$, multiplied by $\gamma$,
are $\operatorname{div}(f\gamma e)$ and $\operatorname{div}(-f^2\langle
e,x\rangle\gamma e)$, since $\partial_e\gamma=-\langle e,x\rangle\gamma$; and
$\int_E\operatorname{div}(h\gamma e)=\delta V_E(h)$ by the Gauss--Green formula.
\end{proof}

For a cluster $U$ and $e\in\R^n$ write $m_i$, $n_i$, $p_i$, $\nabla^\top_i$,
$\ell_i$, $\delta V_i$, $\delta P_i,\dots$ for the objects of
Lemma~\ref{lem:rankone}(b) with $E=U_i$, and
$\delta V(h)=(\delta V_i(h))_i$. As $\sum_i\delta V_i(h)=\int_{\R^n}\operatorname
{div}(h\gamma e)=0$, the vectors $\delta V(h)$ and $\delta^2V(f)$ lie in
$\R^3_0$. A pair $(q,H)$ of $q\in\R^3_0$ and a quadratic form $H$ on $\R^3_0$ is
a \emph{subjet} of $I_w$ at $v\in\mathcal V^\circ$ if
\[
 I_w(v+y)\ge I_w(v)+\langle q,y\rangle+\tfrac12H(y,y)-o(|y|^2)\qquad
 (y\to0,\ y\in\R^3_0).
\]
If $J'$ is $C^2$ near $v$ in the plane of $\mathcal V$ and $I_w-J'$ has a local
minimum at $v$, then $(\nabla J'(v),\nabla^2J'(v))$ is a subjet.

The next lemma is the stability inequality of a minimizer under the variations
$R_t$, with the change of volumes corrected by the subjet. Part (c) admits smooth
test functions that are bounded and constant along the interfaces, but not
compactly supported. Examples are constants, and smooth functions of the common
normal coordinate when the interfaces lie in parallel hyperplanes.

\begin{lemma}\label{lem:touch}
Let $U$ be a minimizer with $v=\gamma(U)\in\mathcal V^\circ$, let $(q,H)$ be a
subjet of $I_w$ at $v$, and $e\in\R^n$.
\begin{enumerate}[label=\textup{(\alph*)}]
\item For $h\in C^\infty_c(\R^n)$,
$\sum_iw_i\,\delta P_i(h)=\langle q,\delta V(h)\rangle$.
\item For $f\in C^\infty_c(\R^n)$,
$H(\delta V(f),\delta V(f))\le\sum_iw_i\int
p_i^2\bigl(|\nabla^\top_if|^2-f^2\bigr)\,dm_i$.
\item Let $u\in C^\infty(\R^n)$ be bounded $m_i$-almost everywhere, with
$\nabla^\top_iu=0$ $m_i$-almost everywhere, for each $i$. Then
\[
 H(y,y)\le-\sum_iw_i\int p_i^2u^2\,dm_i,\qquad y_i=\int u\,p_i\,dm_i .
\]
\end{enumerate}
\end{lemma}

\begin{proof}
By Lemma~\ref{lem:rankone}(a), $R_tU$ is a cluster, so
$I_w(\gamma(R_tU))\le F_w(R_tU)$; and $P_\gamma(U_i)\le F_w(U)/w_i<\infty$. By
Lemma~\ref{lem:rankone}(b), $y(t)=\gamma(R_tU)-v=t\,\delta V(f)+\frac12t^2
\delta^2V(f)+O(t^3)$ and $F_w(R_tU)=F_w(U)+t\sum_iw_i\delta P_i(f)+\frac12t^2
\sum_iw_i\delta^2P_i(f)+O(t^3)$. As $|y(t)|=O(t)$ and $F_w(U)=I_w(v)$, the subjet
inequality gives
\[
 t\Bigl(\langle q,\delta V(f)\rangle-\sum_iw_i\delta P_i(f)\Bigr)
 +\tfrac12t^2\Bigl(\langle q,\delta^2V(f)\rangle+H(\delta V(f),\delta V(f))
 -\sum_iw_i\delta^2P_i(f)\Bigr)\le o(t^2).
\]
Letting $t\to0^\pm$ gives (a) for $h=f$, and then
$H(\delta V(f),\delta V(f))\le\sum_iw_i\delta^2P_i(f)-\langle q,\delta^2V(f)\rangle$.
For $h=f^2\langle e,x\rangle$, $\ell_i(h)=2f\langle e,x\rangle\ell_i(f)
+f^2(|e|^2-p_i^2)$, so $\delta^2P_i(f)+\delta P_i(h)=\int
p_i^2(|\nabla^\top_if|^2-f^2)\,dm_i$, the unbounded terms canceling; and by (a),
$\sum_iw_i\delta P_i(h)=\langle q,\delta V(h)\rangle=-\langle
q,\delta^2V(f)\rangle$. This proves (b). For (c) let $z\in C^\infty_c(\R^n)$,
$0\le z\le1$, $z=1$ on $B_1$, and apply (b) to $f=z(\cdot/R)u$: then
$\nabla^\top_if=u\nabla^\top_i(z(\cdot/R))=O(1/R)$ $m_i$-almost everywhere, and
dominated convergence, with $m_i(\R^n)<\infty$, gives (c) as $R\to\infty$.
\end{proof}

The next lemma recognizes a weighted tripod cluster by its normals. The normals are
constant on each interface and of the form $n_{ij}=(c'_j-c'_i)/w_{ij}$, and the
interfaces have positive measure. The lemma replaces the regularity theory of
minimizers in the equality case of Theorem~\ref{thm:wgdb}.

\begin{lemma}\label{lem:tripod}
Let $U$ be a cluster with $|U_i|>0$ and $a_{ij}>0$ for all $i\ne j$, and let
$c'_1,c'_2,c'_3\in\R^n$ satisfy $c'_j-c'_i=w_{ij}n_{ij}$ $\mathcal
H^{n-1}$-almost everywhere on $S_{ij}$ for all $i\ne j$. Then $U$ is a
weighted tripod cluster up to null sets.
\end{lemma}

\begin{proof}
As $a_{ij}>0$, $|c'_i-c'_j|=w_{ij}$; by the strict triangle inequality the
$c'_i$ are affinely independent.

\emph{Step 1.} Let $c=\sum_lc'_l\one_{U_l}$. For $b\in\R^n$ and $h\in
C^\infty_c(\R^n)$ with $h\ge0$, the Gauss--Green formula, $n_{l'l}=-n_{ll'}$ and
$c'_l-c'_{l'}=-w_{ll'}n_{ll'}$ on $S_{ll'}$ give
\[
 \int\langle c,b\rangle\,\partial_bh\,dx=\sum_l\langle c'_l,b\rangle
 \int_{U_l}\operatorname{div}(hb)\,dx
 =-\sum_{l<l'}w_{ll'}\int_{S_{ll'}}h\langle b,n_{ll'}\rangle^2\,d\mathcal
 H^{n-1}\le0 .
\]
For a mollifier $k_\eps\ge0$ and $c_\eps=c*k_\eps$, taking $h=k_\eps(x-\cdot)$
gives $\partial_b\langle c_\eps,b\rangle(x)\ge0$; integrating along the segment
from $x$ to $y$ with $b=y-x$ gives $\langle c_\eps(y)-c_\eps(x),y-x\rangle\ge0$.
At Lebesgue points $c_\eps\to c$, so $\langle c(y)-c(x),y-x\rangle\ge0$ for
almost every $(x,y)\in\R^n\times\R^n$.

\emph{Step 2.} For $i\ne j$ and almost every $(x,y)\in U_i\times U_j$, Step~1
gives $\langle c'_j-c'_i,x\rangle\le\langle c'_j-c'_i,y\rangle$. For $i<j$ let
$l_{ij}$ be the essential supremum of $\langle c'_j-c'_i,\cdot\rangle$ on $U_i$,
which is finite, and $l_{ji}=-l_{ij}$. Then $\langle c'_j-c'_i,x\rangle\le l_{ij}$
for almost every $x\in U_i$, for all $i\ne j$.

\emph{Step 3.} Let $r_1=\langle c'_2-c'_1,x\rangle-l_{12}$, $r_2=\langle
c'_3-c'_2,x\rangle-l_{23}$, $r_3=\langle c'_1-c'_3,x\rangle-l_{31}$, so that
$r_1+r_2+r_3=-d$ with $d=l_{12}+l_{23}+l_{31}$. By Step~2, up to null sets,
$U_1\subset\{r_1\le0\le r_3\}$, $U_2\subset\{r_2\le0\le r_1\}$ and
$U_3\subset\{r_3\le0\le r_2\}$. The map $x\mapsto(r_1,r_2)$ is an affine surjection
onto $\R^2$. If $d>0$, the nonempty open set $\{r_1,r_2,r_3<0\}$ has positive
measure and meets each cell in a null set. If $d<0$, the same holds for
$\{r_1,r_2,r_3>0\}$. So $d=0$, and $l_{ij}=l_j-l_i$ with $l_1=0$, $l_2=l_{12}$,
$l_3=l_{12}+l_{23}$. By Step~2, $U_i\subset K_i=\{x:\langle c'_i,x\rangle-l_i\ge\langle
c'_j,x\rangle-l_j\ \text{for all }j\}$ up to a null set. The tripod cell $U'_i$,
defined by the strict inequalities for $j\ne i$, misses $K_j$ for $j\ne i$, so
$U'_i\subset U_i$ up to a null set; and $K_i\setminus U'_i$ lies in finitely many
hyperplanes. So $U_i=U'_i$ up to null sets.
\end{proof}

The next lemma computes the interface measures of the sector clusters. It also shows
that the translation $x$ parametrizes the volume vectors in $\mathcal V^\circ$
diffeomorphically, so that the cell volumes determine the translation.

\begin{lemma}\label{lem:V}
For $x,b\in\R^2$ and $i\ne j$,
\begin{gather}
 \tilde a_{ij}(x)=\gamma_1(\langle x,n_{ij}\rangle)\bigl(1-N(\langle
 x,d_{ij}\rangle)\bigr)>0,\label{eq:ray}\\
 \partial_bV_i(x)=\sum_{l\ne i}\langle b,n_{il}\rangle\tilde a_{il}(x)
 =-\bigl(L_{\tilde a(x)/w}Ab\bigr)_i,\label{eq:DV}
\end{gather}
and $V$ is a diffeomorphism of $\R^2$ onto $\mathcal V^\circ$.
\end{lemma}

\begin{proof}
For orthonormal $e,d\in\R^2$ and $x\in\R^2$,
$\int_{x+\R_+d}\gamma_2\,d\mathcal H^1=\gamma_1(\langle
x,e\rangle)(1-N(\langle x,d\rangle))$, since
$|x+rd|^2=\langle x,e\rangle^2+(r+\langle x,d\rangle)^2$; this gives
\eqref{eq:ray}.

The map $V$ is smooth with
values in $\mathcal V^\circ$. If $y\ne x$, then $y-x\in\bar T_i\setminus\{0\}$
for some $i$, so $\langle c_i-c_l,y-x\rangle\ge0$ for $l\ne i$, with strict
inequality for some $l=j$, as $c_i-c_j$ and $c_i-c_k$ are linearly independent.
Then $y+T_i\subset x+T_i$, and the difference contains the nonempty open set of
$z\in x+T_i$ with $\langle c_i-c_j,z-x\rangle<\langle c_i-c_j,y-x\rangle$; so
$V_i(y)<V_i(x)$, and $V$ is injective. For $b\in\R^2$,
$\partial_bV_i(x)=\int_{x+T_i}\operatorname{div}(\gamma_2b)$, and
Lemma~\ref{lem:wedge}(b), applied to $\gamma_2b$ after a cut-off and dominated
convergence, gives \eqref{eq:DV}.
Since the $c_i-c_j$ span $\R^2$, $Ab$ is a multiple of $(1,1,1)$ only if $b=0$.
Moreover $L_{\tilde a(x)/w}$ has kernel $\R(1,1,1)$ on $\R^3$ and is positive
definite on $\R^3_0$. So
$DV(x)\colon\R^2\to\R^3_0$ is an isomorphism, and $V$ is a diffeomorphism onto
an open subset of $\mathcal V^\circ$. Let $\bar c=\frac13\sum c_i$ and
$\hat c_i=(c_i-\bar c)/|c_i-\bar c|$. As $\sum_i|c_i-\bar c|\hat c_i=0$ and the
$\hat c_i$ span $\R^2$, $\max_i\langle\hat c_i,x\rangle\ge m_0|x|$ for some
$m_0>0$, so $K_R=\{x:\langle\hat c_i,x\rangle<R\ \text{for all } i\}$ is a
bounded open triangle containing $0$, on whose boundary
$\max_i\langle\hat c_i,x\rangle=R$. Every $z\in T_i$ has $\langle
c_i-\bar c,z\rangle>0$, so $V_i(x)\le1-N(\langle\hat c_i,x\rangle)$. Let
$\mathcal V_R=\{v\in\mathcal V:\min_iv_i>1-N(R)\}$, a convex set. Then
$V(\partial K_R)\cap\mathcal V_R=\emptyset$, so $\mathcal V_R\cap
V(K_R)=\mathcal V_R\cap V(\bar K_R)$ is open and closed in $\mathcal V_R$; it
contains $V(0)$ for large $R$, hence equals $\mathcal V_R$. As $\mathcal
V^\circ=\bigcup_R\mathcal V_R$, $V$ is onto.
\end{proof}

The \emph{model profile} is $J=\tilde F\circ V^{-1}$ on $\mathcal V^\circ$
(Lemma~\ref{lem:V}).

We now come to the Gaussian double bubble theorem with unequal weights. It says
that the weighted tripod clusters are the only minimizers with volumes in
$\mathcal V^\circ$.

\begin{theorem}[weighted Gaussian double bubble]\label{thm:wgdb}
Let $n\ge2$ and $w_1,w_2,w_3>0$. Then $I_w=J$ on $\mathcal V^\circ$, and every
minimizer $U$ with $\gamma(U)\in\mathcal V^\circ$ is a weighted tripod cluster
up to null sets.
\end{theorem}

\begin{proof}
\emph{Step 1.} Let $U$ be a minimizer with
$v=\gamma(U)\in\mathcal V^\circ$ (Lemma~\ref{lem:exist}(b)), $(q,H)$ a subjet of
$I_w$ at $v$, and $L=L_{a/w}$, positive definite by Lemma~\ref{lem:exist}(c).
Let $\bar n_{ij}=a_{ij}^{-1}\int_{S_{ij}}n_{ij}\,d\gamma^{n-1}$ if
$a_{ij}>0$, $\bar n_{ij}=0$ otherwise, and
$M=\sum_{i<j}a_{ij}(e_i-e_j)\bar n_{ij}^T\colon\R^n\to\R^3_0$. Lemma~\ref{lem:touch}(c)
with $u=1$ has $y=Me$, and its right-hand side is
$-\sum_{i<j}w_{ij}\int_{S_{ij}}\langle e,n_{ij}\rangle^2d\gamma^{n-1}$, as
$w_{ij}=w_i+w_j$. So
\begin{equation}\label{eq:transl}
 H(Me,Me)\le-\sum_{i<j}w_{ij}\int_{S_{ij}}\langle e,n_{ij}\rangle^2\,
 d\gamma^{n-1}\qquad(e\in\R^n).
\end{equation}

\emph{Step 2.} Every $y\in\R^3_0$ satisfies $H(y,y)\le-\langle L^{-1}y,y\rangle$. For
$e\in\R^n$ let $z=L^{-1}Me$ and $N_w=\sum_{i<j}w_{ij}a_{ij}\bar n_{ij}\bar
n_{ij}^T$. Expanding the square and using $\sum_{i<j}a_{ij}\langle\bar
n_{ij},e\rangle(z_i-z_j)=\langle Me,z\rangle$ and $Lz=Me$,
\[
 \sum_{i<j}\frac{a_{ij}}{w_{ij}}\bigl(w_{ij}\langle\bar n_{ij},e\rangle-(z_i-z_j)\bigr)^2
 =e^TN_we-\langle L^{-1}Me,Me\rangle .
\]
Since $\int_{S_{ij}}(n_{ij}-\bar n_{ij})\,d\gamma^{n-1}=0$, the right-hand
side of \eqref{eq:transl} is $-e^TN_we-\sum_{i<j}w_{ij}\int_{S_{ij}}\langle
e,n_{ij}-\bar n_{ij}\rangle^2d\gamma^{n-1}$. Hence
\begin{equation}\label{eq:chain}
 H(Me,Me)+\sum_{i<j}w_{ij}\!\int_{S_{ij}}\!\langle e,n_{ij}-\bar
 n_{ij}\rangle^2d\gamma^{n-1}+\sum_{i<j}\frac{a_{ij}}{w_{ij}}\bigl(w_{ij}\langle\bar
 n_{ij},e\rangle-(z_i-z_j)\bigr)^2\le-\langle L^{-1}Me,Me\rangle,
\end{equation}
which gives the claim for $y$ in the range of $M$. Let $y\notin M(\R^n)$. Then
$M(\R^n)$ and $(\ker M)^\perp$ have dimension at most one. For $e\in\ker M$,
\eqref{eq:transl} gives $\langle e,n_{ij}\rangle=0$ almost everywhere on
$S_{ij}$. Some $a_{ij}$ is positive, so $(\ker M)^\perp=\R\hat e$ for a unit
vector $\hat e$, and $n_{ij}=\eps_{ij}\hat e$ almost everywhere on $S_{ij}$,
with Borel functions $\eps_{ij}=-\eps_{ji}\in\{\pm1\}$. Let $l(x)=\langle
x,\hat e\rangle$.

(i) In Lemma~\ref{lem:touch}(a) with $e=\hat e$, $\ell_i(h)=0$ $m_i$-almost
everywhere, so $-\sum_iw_i\int hl\,dm_i=\sum_iq_i\int h\langle\hat e,n_i\rangle\,
dm_i$ for $h\in C^\infty_c(\R^n)$. Grouped by interfaces, this says that the
measure $\sum_{i<j}(w_{ij}l+(q_i-q_j)\eps_{ij})\,\gamma^{n-1}\llcorner S_{ij}$
vanishes. So $l$ takes values in the finite set
$\mathcal S=\{\pm(q_i-q_j)/w_{ij}\}$ almost everywhere on $\bigcup S_{ij}$,
hence $m_i$-almost everywhere for each $i$.

(ii) Let $s_1<\dots<s_r$ be the elements of $\mathcal S$, $s_0=-\infty$,
$s_{r+1}=\infty$. By (i), no cell has reduced boundary in the connected slab
$O_t=l^{-1}((s_t,s_{t+1}))$, up to a null set; so exactly one cell $U^{(t)}$
has full measure in $O_t$. If $U^{(t-1)}\ne U^{(t)}$, then in
$l^{-1}((s_{t-1},s_{t+1}))$ these two cells are $\{l<s_t\}$ and $\{l>s_t\}$ up to
null sets, and the third cell is null. So, by Lemma~\ref{lem:wedge}(a),
$l^{-1}(s_t)$ lies, up to a null set, in their common interface, on which the
normal pointing from $U^{(t-1)}$ to $U^{(t)}$ is $\hat e$. If $U^{(t-1)}=U^{(t)}$,
no reduced boundary meets $l^{-1}((s_{t-1},s_{t+1}))$. So, up to null sets,
$S_{ij}$ is the union of the hyperplanes $l^{-1}(s_t)$ with
$\{U^{(t-1)},U^{(t)}\}=\{U_i,U_j\}$, and on each the normal pointing from
$U^{(t-1)}$ to $U^{(t)}$ is $\hat e$.

(iii) Let $z=L^{-1}y$ and $b_{ij}=(z_i-z_j)/w_{ij}$, let $f$ be a polynomial
with $f(s_t)=b_{ij}$ whenever $U^{(t-1)}=U_i$ and $U^{(t)}=U_j$ with $i\ne j$, and $u=f\circ l$. Then
$u$ takes finitely many values $m_i$-almost everywhere, $\nabla u$ is parallel to
$\hat e$, and by (ii) $u\langle\hat e,n_{ij}\rangle=b_{ij}$ almost everywhere on
$S_{ij}$, for all $i\ne j$. Lemma~\ref{lem:touch}(c) with $e=\hat e$ has
$\int u\,p_i\,dm_i=\sum_{j\ne i}a_{ij}b_{ij}=(Lz)_i=y_i$ and right-hand side
$-\sum_{i<j}w_{ij}a_{ij}b_{ij}^2=-\langle Lz,z\rangle=-\langle L^{-1}y,y\rangle$.

\emph{Step 3.} By \eqref{eq:ray} and
$\gamma_1'(s)=-s\gamma_1(s)$, $\nabla\tilde a_{ij}(x)=-\langle
x,n_{ij}\rangle\tilde a_{ij}(x)n_{ij}-\gamma_2(x)d_{ij}$, and Young's law gives $\partial_b\tilde F(x)=-\sum_{i<j}w_{ij}\tilde a_{ij}(x)\langle
x,n_{ij}\rangle\langle b,n_{ij}\rangle=-\langle L_{\tilde a(x)/w}Ax,Ab\rangle$. For
$y\in\R^3_0$ and $b=DV(x)^{-1}y$, \eqref{eq:DV} gives $L_{\tilde a(x)/w}Ab=-y$, so
$DJ(V(x))y=\partial_b\tilde F(x)=\langle Ax,y\rangle$: the gradient $\nabla J(V(x))$ is
the orthogonal projection of $Ax$ to $\R^3_0$. Differentiating
$v\mapsto\langle AV^{-1}(v),y\rangle$, and using that $ADV(x)^{-1}y'$ equals
$-L^{-1}_{\tilde a(x)/w}y'$ modulo $\R(1,1,1)$ by \eqref{eq:DV},
\begin{equation}\label{eq:hess}
 \nabla^2J(V(x))=-L^{-1}_{\tilde a(x)/w}\qquad\text{on }\R^3_0 .
\end{equation}

By Lemma~\ref{lem:wedge}(c), $I_w(V(x))\le \tilde F(x)$, so
$I_w\le J$ on $\mathcal V^\circ$. Moreover $J$ extends continuously to $\mathcal
V$ with $J=I_w$ on $\partial\mathcal V$, where $I_w$ is continuous by
Lemma~\ref{lem:exist}(d). Indeed, let $v^h\in\mathcal V^\circ$ converge to
$v\in\partial\mathcal V$, $x^h=V^{-1}(v^h)$, so that $|x^h|\to\infty$ by
Lemma~\ref{lem:V}, and let
$r^h_{ij}$ be the distance from $0$ to the ray $x^h+(\bar T_i\cap\bar T_j)$. If
$z_0$ is the point of a ray nearest to $0$, then $|z|^2\ge|z_0|^2+|z-z_0|^2$ on
the ray, so
\begin{equation}\label{eq:far}
 \tilde a_{ij}(x^h)\le\gamma_1(r^h_{ij}).
\end{equation}
As $r^h_{ij}\ge|\langle x^h,n_{ij}\rangle|$ and the $n_{ij}$ are pairwise
independent, no two of the $r^h_{ij}$ stay bounded along a subsequence. Pass to
a subsequence along which every $r^h_{ij}$ converges in $[0,\infty]$. If all
limits are infinite, then $J(v^h)\to0$ by \eqref{eq:far}, and
$\min_lw_l\sum_iI_\gamma(v^h_i)\le I_w(v^h)\le J(v^h)$ shows that $v$ is a vertex,
where $I_w(v)=0$. Otherwise exactly one, say $r^h_{ij}$, stays bounded. Write
$x^h=s_hn_{ij}+r_hd_{ij}$. Then $r_h\to-\infty$, since $r_h>0$ would give
$r^h_{ij}=|x^h|$, and $s_h\to s$ along a subsequence. A point $z$ with $\langle
z,n_{ij}\rangle<s$ lies in $x^h+T_i$ for large $h$, and one with $\langle
z,n_{ij}\rangle>s$ in $x^h+T_j$: the cells $T_i$ and $T_j$ contain the two halves
of a narrow cone around $d_{ij}$. By dominated convergence $v_i=N(s)$,
$v_j=1-N(s)$, $v_k=0$, while \eqref{eq:ray} and \eqref{eq:far} give
$J(v^h)\to w_{ij}\gamma_1(s)=I_w(v)$.

\emph{Step 4.} Let $\mathcal Q(y)=|\sum_iy_ic_i|^2$ on $\R^3_0$, positive definite
as the $c_i$ are affinely independent, with $\mathcal Q(e_i-e_j)=w_{ij}^2$; for a symmetric
operator $S$ of $\R^3_0$ let $\langle \mathcal Q,S\rangle=\operatorname{tr}(\mathcal QS)$. If
$S=\sum_ls_ly_ly_l^T$, then $\langle \mathcal Q,S\rangle=\sum_ls_l\mathcal Q(y_l)$; so $\langle
\mathcal Q,\cdot\rangle$ is monotone, and $\langle \mathcal Q,S\rangle=0$ with $S\ge0$ forces $S=0$.
Moreover $\langle \mathcal Q,L_b\rangle=\sum_{i<j}w_{ij}^2b_{ij}$, so $\langle
\mathcal Q,L_{a/w}\rangle=F_w(U)$ for every cluster with $F_w(U)<\infty$, and by
\eqref{eq:hess}
\begin{equation}\label{eq:Q}
 \bigl\langle \mathcal Q,(-\nabla^2J(V(x)))^{-1}\bigr\rangle=\langle
 \mathcal Q,L_{\tilde a(x)/w}\rangle=J(V(x)).
\end{equation}
The function $I_w-J$ is lower semicontinuous on $\mathcal V$ and attains its
minimum at some $v_0$. Suppose that it is negative. As $I_w=J$ on
$\partial\mathcal V$, $v_0\in\mathcal V^\circ$. Let $U$ be a minimizer with $\gamma(U)=v_0$. The pair
$(\nabla J(v_0),\nabla^2J(v_0))$ is a subjet of $I_w$ at $v_0$, so Step~2 gives
$-\nabla^2J(v_0)\ge L^{-1}_{a/w}>0$, hence $(-\nabla^2J(v_0))^{-1}\le L_{a/w}$.
Pairing with $\mathcal Q$, \eqref{eq:Q} gives $J(v_0)\le F_w(U)=I_w(v_0)$, a
contradiction. So $I_w\ge J$, and $I_w=J$.

\emph{Step 5.} Let $U$ be a minimizer with $v=\gamma(U)\in\mathcal
V^\circ$, $x_v=V^{-1}(v)$ and $L=L_{a/w}$. As $I_w=J$ is smooth,
$(\nabla J(v),\nabla^2J(v))$ is a subjet at $v$, and Step~2 gives
$D=L-(-\nabla^2J(v))^{-1}\ge0$, with $\langle \mathcal Q,D\rangle=F_w(U)-J(v)=0$ by
\eqref{eq:Q}. So $D=0$ and $L_{a/w}=L_{\tilde a(x_v)/w}$ by \eqref{eq:hess}. Both
matrices vanish on $(1,1,1)$ and have off-diagonal entries $-a_{ij}/w_{ij}$ and
$-\tilde a_{ij}(x_v)/w_{ij}$, so $a_{ij}=\tilde a_{ij}(x_v)>0$. Now \eqref{eq:chain} with
$H=-L^{-1}$ forces both sums on its left to vanish for every $e$. The first gives
$n_{ij}=\bar n_{ij}$ almost everywhere on $S_{ij}$, and the second
$w_{ij}\langle\bar n_{ij},e\rangle=z_i(e)-z_j(e)$, where $z(e)=L^{-1}Me$ is
linear in $e$. Define $c'_i\in\R^n$ by $\langle c'_i,e\rangle=-z_i(e)$. Then
$c'_j-c'_i=w_{ij}n_{ij}$ almost everywhere on $S_{ij}$, and $|U_i|>0$ as
$v\in\mathcal V^\circ$. Lemma~\ref{lem:tripod} completes the proof.
\end{proof}

\begin{remark}
For equal weights Theorem~\ref{thm:wgdb} is the Gaussian double bubble theorem
of Milman and Neeman \cite{MNg}, and the proof follows their strategy. A matrix
inequality for functions touching the profile from below (Step~2) is compared with
the model profile by a maximum principle (Step~4). Unequal weights enter through
the pairing with $\mathcal Q$. No regularity of minimizers is used. The variations
are the rank-one maps $R_t$, and the equality case is settled by
Lemma~\ref{lem:tripod}.
\end{remark}

\begin{remark}\label{rem:lawlor}
In $\R^n$ with the Lebesgue measure, the double bubble problem with unequal
weights describes immiscible fluids. Lawlor \cite{Lawlor} proved that the
standard weighted double bubbles minimize the weighted area among piecewise
smooth boundaries of pairs of open regions with prescribed volumes, and that they
are the only minimizers. If both prescribed volumes are positive and the weights
satisfy the strict triangle inequality, a standard weighted double bubble consists of three caps of $(n-1)$-spheres, one of
which may be a flat disk, with a common boundary $(n-2)$-sphere. Along this
sphere the conormals, with lengths equal to the weights, add up to zero.
Theorem~\ref{thm:wgdb} is the analogue for the Gaussian measure, with the
weighted tripod clusters in place of the standard weighted double bubbles.
\end{remark}

We apply the theorem to the cones over subsets of $S^3$. The next proposition
gives their volumes and Gaussian perimeters.

\begin{proposition}\label{prop:cone}
Let $(X_1,X_2,X_3)$ be a partition of $S^3$, up to null sets, into sets of
finite perimeter, and $K_i=\{rx:r>0,\ x\in X_i\}\subset\R^4$. Then
$K=(K_1,K_2,K_3)$ is a cluster with $\gamma(K_i)=|X_i|/2\pi^2$ and
$P_\gamma(K_i)=b\,P(X_i)$, where $b=(2\pi)^{-2}\int_0^\infty r^2e^{-r^2/2}dr$.
\end{proposition}

\begin{proof}
The volumes follow in polar coordinates $y=rx$, since $\gamma(\R^4)=1$. Let
$Y\in C^1_c(\R^4;\R^4)$ with $|Y|\le1$. Replacing $Y$ by $Y$ times a cut-off
vanishing on $B_\eps$ changes $\int_{K_i}\operatorname{div}(\gamma Y)$ by
$O(\eps^3)$, so let $Y$ vanish near $0$, and write $Y(rx)=f(r,x)x+Z(r,x)$ with
$Z(r,x)\in T_xS^3$. The radial part contributes
$\int_{X_i}\int_0^\infty\partial_r(r^3\gamma f)\,dr\,d\mathrm{vol}=0$,
where $\gamma(r)=(2\pi)^{-2}e^{-r^2/2}$. The tangential part contributes
$\int_0^\infty r^2\gamma(r)\int_{X_i}\operatorname{div}_{S^3}Z(r,\cdot)
\,d\mathrm{vol}\,dr\le b\,P(X_i)$, and fields $Y(rx)=t_l(r)Z(x)$ with
$t_l\in C^\infty_c((0,\infty))$, $0\le t_l\le1$, $t_l\to1$, give values tending
to $b\int_{X_i}\operatorname{div}_{S^3}Z\,d\mathrm{vol}$, whose supremum over
$|Z|\le1$ is $bP(X_i)$. So $P_\gamma(K_i)=bP(X_i)<\infty$,
and $K$ is a cluster.
\end{proof}

The next corollary is the statement on partitions that is used in the proof of
Theorem~D.

\begin{corollary}\label{cor:fan}
Let $a_1,a_2,a_3\in(0,\pi)$ with $a_1+a_2+a_3=2\pi$, and let $w_i$ be
proportional to $\sin a_j+\sin a_k-\sin a_i$, with $w_1+w_2+w_3=1$. Then
$w_i>0$, and every partition $(X_1,X_2,X_3)$ of $S^3$, up to null sets, into
sets of finite perimeter with $|X_i|=\pi a_i$ satisfies
$\sum_iw_iP(X_i)\ge4\pi$, with equality only for the three lunes of a $3$-fan, up
to null sets.
\end{corollary}

\begin{proof}
As $\sin a_i=-\sin(a_j+a_k)$, $\sin a_j+\sin a_k-\sin a_i=\sin a_j(1+\cos a_k)+\sin
a_k(1+\cos a_j)>0$, and $w_{ij}=w_i+w_j$ is proportional to $\sin a_k$. By the
law of sines the triangle $c_1c_2c_3$ has the angle $\pi-a_k$ at $c_k$, so $T_k$
has the angle $a_k$. With $C=\{z_2=0\}\subset\C^2$, the cluster
$(T_i\times\R^2)_i$ is, after an orthogonal map, the cluster of cones over the
three lunes of a $3$-fan of $C$ with angles $a_i$. By Lemma~\ref{lem:wedge}(c)
and Proposition~\ref{prop:cone}, $V(0)=(a_i/2\pi)_i$ and
$\tilde F(0)=b\sum_iw_i\cdot4\pi=4\pi b$. By Proposition~\ref{prop:cone}
and Theorem~\ref{thm:wgdb} with $n=4$, the cones $K$ over the $X_i$ satisfy
\[
 b\sum_iw_iP(X_i)=F_w(K)\ge I_w(V(0))=J(V(0))=4\pi b .
\]
In case of equality $K$ is a minimizer, hence a weighted tripod cluster with
volumes $V(0)$. Up to an orthogonal map it is $((x+T_i)\times\R^2)_i$, and
$x=0$ by injectivity of $V$ (Lemma~\ref{lem:V}). So the $K_i$ are the cones over the three lunes of a
$3$-fan, and the $X_i$ are these lunes up to null sets, since
$|(K_i\triangle K'_i)\cap B_1|=\frac14|X_i\triangle X'_i|$ for cones $K_i,K'_i$
over $X_i,X'_i$.
\end{proof}

\begin{remark}
These numbers $w_i$ are the ones for which the $3$-fan with angles $a_i$ is
stationary for $\sum_iw_iP(X_i)$ without volume constraint. This is Young's law at
$C$ for the weights $w_i+w_j$.
\end{remark}

We now prove Theorem~D.

\begin{proof}[Proof of Theorem~D]
By Lemma~\ref{lem:threecell} it suffices to show that a partition
$(Y_1,Y_2,Y_3)$ of $S^3$ into sets of finite perimeter with $|Y_i|<\pi^2$ and
$P(Y_i)\le4\pi$ consists of the three lunes of a $3$-fan. Then $P(Y_i)=4\pi$ and
$E_i=Y_i$ up to null sets. No $Y_i$ is null, since otherwise the other two would
be complementary with volumes less than $\pi^2$. So $a_i=|Y_i|/\pi\in(0,\pi)$,
$\sum a_i=2\pi$, and with the numbers $w_i$ of Corollary~\ref{cor:fan},
\[
 \max_iP(Y_i)\ge\sum_iw_iP(Y_i)\ge4\pi\ge\max_iP(Y_i).
\]
So equality holds, and Corollary~\ref{cor:fan} gives the $3$-fan.
\end{proof}
\section{The equality case}\label{sec:sharp}

In this section we prove Theorem~B and Corollary~C. The main result is
Theorem~\ref{thm:sharp}. Assume $K_g\le1$, $A(g)=2\pi$, and $\inj(p)\ge\pi$ for all
$p$ in a nonempty open set $V$. Unless $g$ is round, Theorem~\ref{thm:sharp}(a)
produces a stable embedded totally geodesic unit sphere through every point of
$V$. Proposition~\ref{prop:cover} rules this out for $V=S^3$. To find the spheres
we compare $(S^3,g)$ with the round sphere by the exponential map at a point $p$,
and use smooth global modifications of the comparison map $G$
(Lemma~\ref{lem:Geps}). These maps carry almost optimal families to the round
sphere, where Sections~\ref{sec:round} and~\ref{sec:threeset} apply. The equality
case of Rauch's theorem and the Codazzi equation then give the sphere.

In this section $g_0$ is the round metric of curvature one. For a closed
two-sided surface with unit normal $N$, $R_N(X,Y)=\langle R(X,N)N,Y\rangle$ for
tangent $X,Y$. A closed totally geodesic surface is \emph{stable} if its Jacobi
operator $-\Delta-\Ric(N,N)$ has no negative eigenvalue. The main result is the
following theorem.

\begin{theorem}\label{thm:sharp}
Let $g$ be a metric on $S^3$ with $K_g\le1$ and $A(g)=2\pi$.
\begin{enumerate}[label=\textup{(\alph*)}]
\item Suppose $\inj(p)\ge\pi$ for all $p$ in a nonempty open set $V$. Then either
$g$ is round, or every point $p$ of $V$ lies on a stable embedded totally
geodesic sphere isometric to the unit sphere, namely the image under $\exp_p$ of
the closed disk of radius $\pi$ in a plane of $T_pS^3$.
\item In particular, if in \textup{(a)} $(S^3,g)$ contains no stable embedded
minimal sphere, for instance if $\Ric_g>0$, then $g$ is round.
\item If $\inj(p)\ge\pi$ for all $p\in S^3$, then $g$ is round.
\end{enumerate}
\end{theorem}

The last step uses the following theorem of Mazet \cite[Corollary~1]{Mazet}.
\emph{A complete Riemannian $n$-manifold, $n\ge3$, with $K\le1$ that admits a
totally geodesic isometric immersion of the unit $(n-1)$-sphere which is
unstable as a minimal hypersurface has universal cover isometric to the unit
$n$-sphere.} The word complete is not in the statement of \cite{Mazet}, but it
is understood: a thin tubular neighbourhood of an equator of the unit sphere
satisfies the other hypotheses. Here only $S^3$ occurs.

Let $g$ be a metric on $S^3$ with $K_g\le1$. Fix $p$ with $\inj(p)\ge\pi$, a
point $q$ of the round sphere and a linear isometry $I\colon T_pS^3\to T_qS^3$.
Let $U=\exp_p(B_\pi)$, the open ball of radius $\pi$ about $p$, and
\[
 G=\exp_q\circ I\circ\exp_p^{-1}\colon U\to S^3\setminus\{-q\}.
\]
Write $y=\exp_q(tv)$, $0<t<\pi$, $|v|=1$, and let $c_v$ be the geodesic of $g$
with initial velocity $I^{-1}v$. Then $dG^{-1}_y$ maps the radial unit vector to
$c_v'(t)$, and maps the parallel transport $P_tw$ of $w\perp v$ along
$s\mapsto\exp_q(sv)$ to $J(t)/\sin t$, where $J$ is the Jacobi field along
$c_v$ with $J(0)=0$, $J'(0)=I^{-1}w$. By the Gauss lemma the two images are
orthogonal.

The next lemma is Rauch's theorem with its equality case.

\begin{lemma}[Rauch, with equality]\label{lem:rauch}
Let $K_g\le1$, let $c\colon[0,\pi]\to S^3$ be a geodesic with unit speed, and
$J$ a Jacobi field along $c$ with $J(0)=0$ and $J'(0)\perp c'(0)$. Then
$|J(t)|\ge\sin t\,|J'(0)|$ on $[0,\pi]$. If equality holds at some
$t_0\in(0,\pi)$, then on $[0,t_0]$ one has $J=\sin t\,e$ with $e$ parallel and
$R(e,c')c'=e$.
\end{lemma}
\begin{proof}
If $J'(0)=0$, then $J=0$, and the claim holds with $e=0$.
Where $J\ne0$, $|J|''\ge-K(c',J)|J|\ge-|J|$, with equality only if $J'$ is
parallel to $J$ and $K(c',J)=1$. So $|J|/\sin t$ is nondecreasing, with limit
$|J'(0)|$ at zero, and equality at $t_0$ gives equality on $[0,t_0]$. Then
$J/|J|$ is parallel, and an eigenvector of $R(\cdot,c')c'$ for the largest
possible eigenvalue.
\end{proof}

So $G$ does not increase lengths. For a plane $P\subset T_yS^3$ let
$j(y,P)\ge1$ be the factor by which $dG^{-1}_y$ multiplies areas in $P$. Let
$l_1,l_2\ge1$ be the singular values of the restriction $T$ of $dG^{-1}_y$ to
the vectors orthogonal to the radial direction. For an orthonormal basis
$X=\cos a\,\partial_r+\sin a\,w_1$, $Y=w_2$ of $P$ with $w_1,w_2\perp\partial_r$,
\begin{equation}\label{eq:j}
 j(y,P)^2=\cos^2\!a\,|Tw_2|^2+\sin^2\!a\,(l_1l_2)^2 .
\end{equation}
For a unit vector $n\in T_yS^3$, $j(y,n^\perp)$ is the length of
$\operatorname{cof}(dG^{-1}_y)\,n$; so $v\mapsto|\operatorname{cof}(dG^{-1}_y)\,v|$
is a norm on $T_yS^3$ that extends $n\mapsto j(y,n^\perp)$.

By \eqref{eq:j}, $dG^{-1}$ does not decrease the area of any plane. The next lemma
is the equality case. Preservation of area on one plane gives equality in
Lemma~\ref{lem:rauch}. On a plane transversal to the radial direction it gives
equality for all Jacobi fields, hence the round metric along the segment from $q$
to $y$.

\begin{lemma}\label{lem:j}
Let $K_g\le1$, $\inj(p)\ge\pi$ and $G$ as above. Suppose $j(y,P)=1$,
$y=\exp_q(tv)$, $0<t<\pi$, let $w_2$ be a unit vector of $P$ orthogonal to the
radial direction, and $P_t$ the parallel transport along $s\mapsto\exp_q(sv)$.
\begin{enumerate}[label=\textup{(\alph*)}]
\item If $P$ does not contain the radial direction, then $(G^{-1})^*g=g_0$ at
every point $\exp_q(sv)$, $0\le s\le t$.
\item If $P$ contains the radial direction, then the Jacobi field along $c_v$
with $J(0)=0$ and $J'(0)=I^{-1}P_t^{-1}w_2$ has equality in
Lemma~\ref{lem:rauch} on $[0,t]$.
\end{enumerate}
\end{lemma}
\begin{proof}
In (a), $\sin a\ne0$ forces $l_1=l_2=1$, so every Jacobi field of
Lemma~\ref{lem:rauch} has equality on $[0,t]$. In polar coordinates
$(G^{-1})^*g$ is determined by their inner products, which are those of the
round sphere. In (b), $\cos a=1$ and $|Tw_2|=1$.
\end{proof}

To transport whole families of loops, we use smooth global modifications of the
comparison map $G$. They are maps of $S^3$ of degree $\pm1$.

\begin{lemma}\label{lem:Geps}
Let $K_g\le1$, $\inj(p)\ge\pi$ and $G$ as above. There is $\eps_0>0$ such that
for every $\eps\in(0,\eps_0)$ there is a smooth map $G_\eps\colon S^3\to S^3$
of degree $\pm1$, $(1+\eps)$-Lipschitz from $g$ to $g_0$, with
$G_\eps=F_\eps\circ G$ on $U$ and $G_\eps=-q$ near $S^3\setminus U$, where
$F_\eps(\exp_q(tv))=\exp_q(f_\eps(t)v)$ and $F_\eps\to\mathrm{id}$ in $C^1$ on
compact subsets of $S^3\setminus\{-q\}$ as $\eps\to0$.
\end{lemma}
\begin{proof}
Take $f=f_\eps$ smooth and nondecreasing, $f(t)=(1+\eps)t$ until $f=\pi-2\eps$,
then $0<f'\le1+\eps$ until $f$ reaches $\pi$ at some $t_2<\pi$, flat to
infinite order there, and $f=\pi$ afterwards. The radial stretch is at most
$1+\eps$. A tangential vector is multiplied by
$\sin f(t)/|J(t)|\le\sin f(t)/\sin t$, which is at most $1+\eps$: on the linear
part this is elementary, and near $t_2$, for small $\eps$,
$\sin f\le\sin((1+\eps)(t_2-t))\le(1+\eps)\sin(\pi-t)$. Near $\{t\ge t_2\}$ the
map is, in normal coordinates at $-q$, a smooth unit vector field times a
function vanishing to infinite order, so $G_\eps$ is smooth. It maps
$\{t<t_2\}$ diffeomorphically onto $S^3\setminus\{-q\}$.
\end{proof}

Unlike $G$, the map $G_\eps$ is defined on all of $S^3$. It has degree $\pm1$ and
increases lengths at most by the factor $1+\eps$. Composing loop families, disks
and sweepouts with $G_\eps$ gives $A(g_0)\le(1+\eps)^2A(g)$ and
$W(g_0)\le(1+\eps)^2W(g)$. If $A(g)=2\pi$, $A(\Gamma_n,g)\to2\pi$ for
families $\Gamma_n$ of degree one and $\eps_n\to0$, then
$A(G_{\eps_n}\circ\Gamma_n,g_0)\to2\pi$.

We now prove Theorem~\ref{thm:sharp}.

\begin{proof}[Proof of Theorem~\ref{thm:sharp}\textup{(a)} and \textup{(b)}]
Let $p\in V$. Take families $\Gamma_n$ of degree one with
$A(\Gamma_n,g)<2\pi+\frac1n$, numbers $\eps_n\in(0,\eps_0)$ with $\eps_n\to0$,
and $G_n=G_{\eps_n}$. Let
$L_n$ bound the lengths of the loops of $\Gamma_n$, and $d_n>0$ with $d_n\to0$
and $d_nL_n\to0$.

\emph{Step 1.} We apply Theorem~\ref{thm:switching} to the loop families
$G_n\circ\Gamma_n$. They have degree $\pm1$, and for each $n$ their loops have a
common Lipschitz bound. We use $\eps=\eps_n$ and a triangulation so fine that the
loops of $\Gamma_n$ at the vertices of a common triangle are at uniform
$g$-distance at most $d_n$. The disks are $G_n\circ f_x$, where $f_x$ is a Lipschitz
disk in $(S^3,g)$ spanning $\Gamma_n(x)$ with $\Area_g(f_x)\le2\pi+\frac2n$,
constant if $\Gamma_n(x)$ is constant. This gives three vertices. Let
$f_1^n,f_2^n,f_3^n$ be their disks in $(S^3,g)$. For each pair
$k\in\{12,23,31\}$ let $s_k^n$ be the switching sphere, $u_k^n$ the bounding
function and $E_k^n$ the odd set of Theorem~\ref{thm:switching}. Then $|E_k^n|\le\pi^2$, the sum of the
$\one_{E_k^n}$ is odd outside a set of volume at most $\eps_n$, and
$\Area_{g_0}(s_k^n)\le2(1+\eps_n)^2(2\pi+\frac2n)+\eps_n$.

\emph{Step 2.} We pass to the limit. By Lemma~\ref{lem:bounding}(d),
$|Du_k^n|(S^3)\le\Area_{g_0}(s_k^n)$. By Proposition~\ref{prop:limits}, after
passing to a subsequence and adding even integers, $u_k^n\to u_k$ in $L^1$ and
$Du_k^n\to Du_k$ weakly-*. The odd sets $E_k$ of $u_k$ satisfy
$P(E_k)\le|Du_k|(S^3)\le4\pi$ and $|E_k|\le\pi^2$, and
$\one_{E_{12}}+\one_{E_{23}}+\one_{E_{31}}$ is odd almost everywhere. In the \emph{lune case} all three volumes are less than $\pi^2$, and
the $E_k$ are, up to null sets, the three closed lunes of a $3$-fan by Theorem~D,
with angles less than $\pi$. In the \emph{hemisphere case} some $|E_k|=\pi^2$,
and $E_k$ is a closed hemisphere up to a null set by Schmidt's inequality
(Section~\ref{sec:round}).
In either case fix $k$ and a closed lune or closed hemisphere $E$ with $E=E_k$
up to a null set. By Proposition~\ref{prop:lune}, $P(E)=4\pi$,
$|Du_k|=\mathcal H^2\llcorner\partial E$ and $u_k=b\pm\one_E$ almost everywhere
for an integer $b$.

\emph{Step 3.} We now show that the area factor is one along $\partial E$, because
no area is lost in the limit. More precisely, $j(y,T_y\partial E)=1$ at every point
$y\ne-q$ of $\partial E$ where $\partial E$ is smooth, that is, off the great
circle of the $3$-fan in the lune case. First,
\[
 4\pi=|Du_k|(S^3)\le\liminf|Du_k^n|(S^3)
 \le\limsup\Area_{g_0}(s_k^n)\le4\pi ,
\]
so the total variations converge. Let $h$ be continuous with $0\le h\le1$ and
compact support in $S^3\setminus\{-q\}$, and let $Z\subset S^3\setminus\{-q\}$
be a compact set containing $\operatorname{supp}h$ and
$F_{\eps_n}^{-1}(\operatorname{supp}h)$ for all large $n$. By Reshetnyak's
continuity theorem \cite{Reshetnyak}, \cite[Theorem~2.39]{AFP}, applied to the
positively homogeneous integrand $h(y)|\operatorname{cof}(dG^{-1}_y)v|$,
\[
 \int h\,j(y,n_{u_k^n}^\perp)\,d|Du_k^n|\longrightarrow
 \int_{\partial E}h\,j(y,T_y\partial E)\,d\mathcal H^2 .
\]
By Lemma~\ref{lem:bounding}(d), applied to this integrand extended by zero to
$-q$, the left side is at most the integral of $h\,j$ over $s_k^n$ counted with
multiplicity. The sphere $s_k^n$ is glued from the disks $G_n\circ f_i^n$,
$G_n\circ f_j^n$ and the annulus of Theorem~\ref{thm:switching}. The annulus has
area at most $\eps_n$, and $j$ is bounded on $\operatorname{supp}h$. Where
$G_n\circ f_i^n\in\operatorname{supp}h$, the disk $f_i^n$ takes values in $U$,
since $G_n=-q$ near $S^3\setminus U$. So $G_n\circ f_i^n=F_{\eps_n}\circ G\circ
f_i^n$ and $G\circ f_i^n\in Z$ there. Hence $h\,j$ times the round Jacobian of
$G_n\circ f_i^n$ is at most the $g$-Jacobian of $(G^{-1}\circ F_{\eps_n}\circ
G)\circ f_i^n$. Moreover $G^{-1}\circ F_{\eps_n}\circ G\to\mathrm{id}$ in $C^1$ on
$G^{-1}(Z)$. The same holds for $f_j^n$. So the integral is at most
$(1+o(1))(\Area_g(f_i^n)+\Area_g(f_j^n))+o(1)\le4\pi+o(1)$. Thus
$\int_{\partial E}hj\le4\pi=\mathcal H^2(\partial E)$ for all such $h$. Letting
$h$ increase to the indicator of $S^3\setminus\{-q\}$ and using that $j\ge1$ is
continuous gives $j=1$ on the smooth part of $\partial E$ away from $-q$.

\emph{Step 4.} We use the equality $j=1$ to find a totally geodesic unit sphere in
$(S^3,g)$. In the lune case Step~3 applies to all three lunes, so $j=1$ on
the three pages away from $-q$. Let $C$ be the great circle of the $3$-fan, $P_1,P_2,P_3$ the
pages and $S_i$ the great sphere containing $P_i$. Since the angles are less
than $\pi$, the half of $S_i$ opposite to $P_i$ lies in the lune not adjacent to
$P_i$. Any two of the $S_i$ meet only in $C$.

\emph{Case 1: $q$ is not on $\partial E$ in the hemisphere case, or $q$ is on
none of the $S_i$, or $q$ lies on a page $P_i$ but not on $C$.} Then
$(G^{-1})^*g=g_0$ on a closed hemisphere $\overline H$ that contains $q$ and not
$-q$. If $E$ is a hemisphere, every ray from $q$ crosses the great sphere
$\partial E$ transversally before reaching $-q$, and Lemma~\ref{lem:j}(a)
applies with $\overline H$ the side of $\partial E$ containing $q$. If $q$ is on
none of the $S_i$, then $q$ and $-q$ are interior to two different lunes. Along
a ray from $q$ avoiding $C$ the angle about $C$ is monotone. So the ray crosses
the pages transversally, the last time when it enters the lune $L'$ containing
$-q$. Lemma~\ref{lem:j}(a) at the last crossing and continuity give
$(G^{-1})^*g=g_0$ on $S^3\setminus\operatorname{int}L'$. This set is the union of
the other two lunes, a lune of angle larger than $\pi$, and it contains a closed
hemisphere not containing $-q$. If $q\in P_i\setminus C$, the rays from $q$ not
contained in $S_i$ cross one of the other pages transversally before reaching $-q$
in the lune not adjacent to $P_i$. Again $(G^{-1})^*g=g_0$ on the union of the
other two lunes, by continuity. So $G^{-1}(\overline H)$ is a round closed hemisphere
in $(S^3,g)$, and its boundary $\Sigma$ is a totally geodesic unit sphere.

\emph{Case 2: $q\in S=\partial E$ in the hemisphere case, or $q\in S=S_i$ on
the half opposite to $P_i$, not on $C$.} Every ray from $q$ inside $S$ runs, on
a final stretch ending at $-q$, in $\partial E$ or in the page $P_i$. By
Lemmas~\ref{lem:j}(b) and~\ref{lem:rauch}, for orthonormal $v,w\in T_qS$ the
Jacobi field along $c_v$ with $J(0)=0$ and $J'(0)=I^{-1}w$ is $\sin t\,e_2$ on
$[0,\pi)$, with $e_2$ parallel and $R(e_2,c_v')c_v'=e_2$. So $G^{-1}$ restricts to an
isometric embedding of $S\setminus\{-q\}$ onto a surface $D$ whose unit normal
$N$ is parallel along the $c_v$. In the frame $e_1=c_v'$, $e_2$ the second
fundamental form of $D$ is $\ell\,e_2^*\otimes e_2^*$, and the Codazzi equation
reads $(\ell\sin t)'=\pm\sin t\,\langle R(e_1,e_2)e_2,N\rangle=0$. Indeed, since
$K\le1=K(e_1,e_2)$, $e_1$ is an eigenvector of $R(\cdot,e_2)e_2$. As $\ell$ is
bounded near $p$, $\ell=0$, and $D$ is totally geodesic. Hence, for $a\in S$ close
to $-q$, $G^{-1}\circ\exp^S_a=\exp_{G^{-1}(a)}\circ dG^{-1}_a$ on the vectors of
$T_aS$ of length less than $\pi$ whose geodesic segment avoids $-q$. The right
side is smooth on the disk of radius $\pi$ in $T_aS$. So $G^{-1}$ extends from $S\setminus\{-q\}$ to a smooth totally geodesic
isometric immersion of $S$, which maps $-q$ to $p^*=\exp_p(\pi I^{-1}v)$ for
every unit $v\in T_qS$. Since $\inj(p)\ge\pi$, $d(p,p^*)=\pi$, while the points
of $D$ have distance less than $\pi$ from $p$. So $\Sigma=D\cup\{p^*\}$ is an
embedded totally geodesic unit sphere through $p$. In the hemisphere case it is
the image under $\exp_p$ of the closed disk of radius $\pi$ in the plane
$I^{-1}T_qS$.

\emph{Case 3: $q\in C$.} Then $E$ is one of the three lunes $E_{12},E_{23},E_{31}$
of a $3$-fan with great circle $C$, all of volume less than $\pi^2$. Choose $m>0$
with $|E_k|<\pi^2-3m$ for each $k$. For $x\in V$ near $p$ let
$G_x=\exp_q\circ I_x\circ\exp_x^{-1}$ on $U_x=\exp_x(B_\pi)$, with linear
isometries $I_x\colon T_xS^3\to T_qS^3$ depending continuously on $x$ and
$I_p=I$. Let $G_{x,\eps}$ be the maps of Lemma~\ref{lem:Geps}. Let $\operatorname{Jac}_x$ and
$\operatorname{Jac}_{x,\eps}$ be the Jacobians of $G_x$ and $G_{x,\eps}$ from $g$ to $g_0$, with
$\operatorname{Jac}_x=0$ outside $U_x$. Since $0\le \operatorname{Jac}_x\le1$, $\exp$ depends smoothly on its base
point and the distance spheres $\partial U_x$ are null sets, $\operatorname{Jac}_x\to \operatorname{Jac}_p$ in
$L^1(S^3,g)$ as $x\to p$. Also $\operatorname{Jac}_{x,\eps}\to \operatorname{Jac}_x$ in $L^1$ as $\eps\to0$, since
$F_\eps$ is $(1+\eps)$-Lipschitz and converges to the identity in $C^1$ away
from $-q$. Choose $p_1\in V$ with $\|\operatorname{Jac}_{p_1}-\operatorname{Jac}_p\|_{L^1}<m$, not on
$G_p^{-1}(C)$, and so close to $p$ that an open ball $B$ about $p$ containing
$p_1$ has compact closure in $U_p\cap U_{p_1}$.

For each edge $k$ of the triangle of Step~1 let $\tilde
a_k^n$ be the annulus in $(S^3,g)$ swept out by the minimizing $g$-geodesics
between corresponding points of the loops of its two vertices. Its $g$-area is
at most a constant times $d_nL_n$, which tends to zero. Let $\tilde s_k^n$ be the Lipschitz sphere glued from the
two disks $f_i^n,f_j^n$ and $\tilde a_k^n$, oriented as the switching sphere
$s_k^n$. Bounding functions and the class $BV(S^3;\Z)$ do not depend on the
metric. Let $v_k^n$ be bounding functions of the $\tilde s_k^n$
(Lemma~\ref{lem:bounding}(c)). The disks cancel in $\sum_k\tilde
s^n_{k\#}[\![S^2]\!]$, so $\sum_kv_k^n$ has total $g$-variation at most
$\sum_k\Area_g(\tilde a_k^n)\to0$. By the isoperimetric inequality of $g$, it
differs from an integer $c_n$ only on a set of volume tending to zero. The map
$G_{x,\eps}$ maps a subset of $U_x$ diffeomorphically onto $S^3\setminus\{-q\}$
and the rest of $S^3$ to $-q$. By the area formula, as in
Lemma~\ref{lem:bounding}(b), the functions $\pm v_k^n\circ G_{x,\eps_n}^{-1}$ are
bounding functions of the spheres $G_{x,\eps_n}\circ\tilde s_k^n$. Here the inverse
is that of the restriction, and the sign is the same for all $k$. Their sum is
$\pm c_n$ outside a set of volume tending to zero. For $x=p$ these spheres
differ from the switching spheres $s_k^n$ only in the annuli, by cycles of
round mass tending to zero. By \eqref{eq:schmidt} and the coarea formula, after
adding integers to the $v_k^n$ and changing $c_n$ accordingly, $u_k^n$ and
$\pm v_k^n\circ G_n^{-1}$ agree outside a set of volume tending to zero. As
$\sum_ku_k^n$ is odd outside a set of volume at most $\eps_n$, $c_n$ is odd for
large $n$.

The odd sets of $\pm v_k^n\circ G_{x,\eps_n}^{-1}$
have volume $\int_{\{v_k^n\ \text{odd}\}}\operatorname{Jac}_{x,\eps_n}$. For $x=p$ these volumes
differ from $|E_k^n|$ by $o(1)$, so they tend to $|E_k|<\pi^2-3m$. Since
$\|\operatorname{Jac}_{p_1}-\operatorname{Jac}_p\|_{L^1}<m$, for $x=p_1$ they are less than $\pi^2-m$ for large
$n$. So the spheres $G_{p_1,\eps_n}\circ\tilde s_k^n$, made of the disks
$G_{p_1,\eps_n}\circ f_i^n$ and the annuli $G_{p_1,\eps_n}\circ\tilde a_k^n$ of
area tending to zero, with the bounding functions $\pm v_k^n\circ
G_{p_1,\eps_n}^{-1}$, satisfy the conclusions of Step~1 with some $\eps_n'\to0$
in place of $\eps_n$. Hence Steps~2 and~3 apply at $p_1$. The limits of the odd
sets have volume at most $\pi^2-m$, so they are the three lunes $E_k'$ of a
$3$-fan with a great circle $C'$. The area factor of $dG_{p_1}^{-1}$ is one on its
pages away from $-q$. The $BV$ norms of $g$ and $g_0$ are equivalent, so the
$L^1$ bound in the proof of Proposition~\ref{prop:limits}, applied to the
$v_k^n$ for $g_0$, gives, after passing to a subsequence and adding even
integers, $v_k^n\to v_k$ in $L^1$. As $G_{x,\eps_n}\to G_x$ in $C^1$ on $\bar B$
for $x=p,p_1$, and composition with diffeomorphisms converging in $C^1$
preserves $L^1$ convergence of indicator functions, the odd set of $v_k$
coincides on $B$, up to null sets, with $G_p^{-1}(E_k)$ and with
$G_{p_1}^{-1}(E_k')$. The points of $B$ every neighbourhood of which meets each
of these three sets in positive volume form $G_p^{-1}(C)\cap B$, and also
$G_{p_1}^{-1}(C')\cap B$. As $p_1\notin G_p^{-1}(C)$, the point
$q=G_{p_1}(p_1)$ is not on $C'$, and the lune subcase of Case~1 or of Case~2
applies at $p_1$.

\emph{Step 5.} It remains to conclude. In every case $(S^3,g)$ contains an embedded
totally geodesic unit sphere $\Sigma$. The metric $g$ is round by Mazet's theorem if
$\Sigma$ is unstable. In Case~1 the sphere bounds a round hemisphere, so $\Ric(N,N)=2$ along
it, and the constants are an unstable direction. In the lune subcase of Case~2,
$q$ is interior to the lune $L$ opposite to $P_i$. The rays from $q$ outside $S$
cross a page of $L$ transversally. So Lemma~\ref{lem:j}(a) gives
$(G^{-1})^*g=g_0$ on $L$, which contains the half of $S$ opposite to $P_i$. So
$\Ric(N,N)=2$ on a closed hemisphere of $\Sigma$. If $\Sigma$ were stable, the
first Dirichlet eigenfunction of that hemisphere, extended by zero, would have
nonpositive second variation. So it would be a first eigenfunction of the Jacobi
operator, which cannot vanish on an open set. Since Case~3 leads to a
lune subcase at $p_1$, $g$ is round unless the hemisphere subcase of Case~2
occurs at $p$ with $\Sigma$ stable, which is the second alternative of (a). If
$\Ric_g>0$, every closed two-sided minimal surface is unstable; this gives (b).
\end{proof}

\begin{remark}
If $q\in C$, every ray from $q$ stays in $C$, in one page, or, apart from its
end points, in the interior of one lune. No ray crosses a page transversally,
so Lemma~\ref{lem:j}(a) does not apply. The source spheres of Case~3, fixed once
the triangle of Step~1 is chosen, do not depend on the base point. This allows the
limits at $p$ and at $p_1$ to be compared.
\end{remark}

If the second alternative of Theorem~\ref{thm:sharp}(a) held for $V=S^3$, the stable
totally geodesic unit spheres given by it would cover $S^3$. The next proposition
excludes this.

\begin{proposition}\label{prop:cover}
Let $g$ be a metric on $S^3$ with $K_g\le1$ that is not round, and $\mathcal S$
the set of embedded totally geodesic spheres in $(S^3,g)$ isometric to the
unit sphere.
\begin{enumerate}[label=\textup{(\alph*)}]
\item Every member of $\mathcal S$ is stable, and every sequence in
$\mathcal S$ has a subsequence that converges smoothly to a member.
\item Two different members are disjoint, or they meet transversally in one
closed geodesic, a great circle of both.
\item A member $\Sigma$ that is a limit of other members is disjoint from every
other member, and there is $w>0$ on $\Sigma$ with $\nabla^2w+wR_N=0$.
\item The members of $\mathcal S$ do not cover $S^3$. If they cover an open set
$V$, then $\overline V$ is covered by pairwise disjoint members as in
\textup{(c)}.
\end{enumerate}
\end{proposition}
\begin{proof}
(a) Stability is Mazet's theorem. Let $h_j\colon S^2\to S^3$ be isometric
embeddings onto members. Since they are totally geodesic,
$h_j\circ\exp_a=\exp_{h_j(a)}\circ dh_j(a)$. A subsequence of $dh_j(a_0)$
converges, so $h_j$ converges smoothly on $S^2\setminus\{-a_0\}$, and then,
using a second point, on $S^2$, to a totally geodesic isometric immersion $h$.
If $h(a)=h(b)$ with $a\ne b$, the two sheets cannot cross transversally,
because the sheets of $h_j$ would. So they are tangent, and the isometry $r$ of
$S^2$ with $r(a)=b$ and $dh_b\,dr_a=dh_a$ satisfies $h\circ r=h$. Such
isometries form a group acting freely, hence of order two, generated by an
orientation-reversing involution, and $h$ would embed the projective plane in
$S^3$. So $h$ is an embedding.

(b) Totally geodesic surfaces tangent at a point coincide, so different
members meet transversally. If two of them meet, the intersection is a union of
disjoint closed geodesics, each a great circle of both. Two great circles of a
round sphere meet, so the intersection is a single one.

(c) Let $\Sigma_j\to\Sigma$, $\Sigma_j\ne\Sigma$. For large $j$, $\Sigma_j$ is
the graph $\{\exp_y(u_jN)\}$ of $u_j\not\equiv0$, $u_j\to0$ smoothly. The second
fundamental form of the graph of $u$ is a second order expression with
linearization $\nabla^2u+uR_N$ at zero. By Schauder estimates a subsequence of
$u_j/\max|u_j|$ converges in $C^2$ to $w\ne0$ with $\nabla^2w+wR_N=0$, so $w$
is a Jacobi field. As $\Sigma$ is stable, $w$ is a first eigenfunction and has
no zeros, so $\Sigma_j\cap\Sigma=\emptyset$ for large $j$. Replace $w$ by $-w$
if necessary. If $\Sigma'\ne\Sigma$ met $\Sigma$, by (b) it would cross $\Sigma$
and, for large $j$, $\Sigma_j$, in great circles of $\Sigma'$, which meet. This
contradicts $\Sigma\cap\Sigma_j=\emptyset$.

(d) With the Hausdorff distance $\mathcal S$ is compact, by (a). The set
$\mathcal S'$ of members that are limits of other members is closed, and its
members are pairwise disjoint by (c). The other members are isolated, hence
countably many, and the complement of a countable union of surfaces is dense. So if
$\mathcal S$ covers an open set $V$, the closed union of $\mathcal S'$ contains
$\overline V$. Suppose $V=S^3$. Fix $\Sigma_0\in\mathcal S'$. Every other
$\Sigma\in\mathcal S'$ bounds a domain $D(\Sigma)$ not containing $\Sigma_0$, of
volume $v(\Sigma)$. Let $\Sigma_j$ have $v(\Sigma_j)\to v_*=\inf v$ and
$\Sigma_j\to\Sigma_*$. If $\Sigma_*=\Sigma_0$, a subsequence of the $\Sigma_j$
consists of graphs over $\Sigma_0$ on one side $B$, so $v(\Sigma_j)<|B|$ and
$v(\Sigma_j)\to|B|=v_*$, against the definition of $v_*$. If
$\Sigma_*\ne\Sigma_0$, then $v(\Sigma_*)=v_*$, and a point of $D(\Sigma_*)$ lies
on some $\Sigma\in\mathcal S'$ disjoint from $\Sigma_0$ and $\Sigma_*$. So
$D(\Sigma)\subset D(\Sigma_*)$ misses a neighbourhood of $\Sigma_*$, and
$v(\Sigma)<v_*$. Both are contradictions.
\end{proof}

\begin{proof}[Proof of Theorem~\ref{thm:sharp}\textup{(c)}]
If $g$ is not round, part (a) with $V=S^3$ says that $\mathcal S$ covers $S^3$,
against Proposition~\ref{prop:cover}(d).
\end{proof}

\begin{proof}[Proof of Theorem~B]
This is Theorem~\ref{thm:sharp}(c).
\end{proof}

\begin{proof}[Proof of Corollary~C]
On a compact simply connected manifold of odd dimension with
$\frac14\le K\le1$ the injectivity radius is at least $\pi$
\cite{CheegerGromoll,KlingenbergSakai}. Apply Theorem~B.
\end{proof}

The next lemma relates the hypothesis of Theorem~B to the lengths of closed
geodesics. It gives Remark~\ref{rem:inj} and is used in
Section~\ref{sec:conj}.

\begin{lemma}\label{lem:kling}
Let $g$ be a metric on $S^3$ with $K_g\le1$.
\begin{enumerate}[label=\textup{(\alph*)}]
\item We have $\inj(g)\ge\pi$ if and only if every closed geodesic has length at
least $2\pi$.
\item If $\inj(g)<\pi$, then a shortest closed geodesic has length
$2\inj(g)$.
\end{enumerate}
\end{lemma}
\begin{proof}
By Klingenberg's lemma, the injectivity radius of a compact manifold with
$K\le1$ is at least the smaller of $\pi$ and half the length of a shortest
nonconstant closed geodesic \cite[Lemma~5.6 and the proof of
Corollary~5.7]{CheegerEbin}. The positive lower curvature bound assumed in
\cite[Corollary~5.7]{CheegerEbin} is used there only through compactness.
Conversely, a closed geodesic of length $\ell$ through $y$ is a geodesic loop at
$y$, so $\ell\ge2\inj(y)\ge2\inj(g)$. These two facts give (a) and (b).
\end{proof}

\begin{remark}\label{rem:MN}
If $\delta\le K_g\le1$ with $0<\delta\le1$, the metric $\delta g$ has positive
Ricci curvature and scalar curvature at least $6$. By Marques and Neves,
$\widetilde W(\delta g)\le4\pi$, with equality only if $\delta g$ has constant
sectional curvature one \cite[Theorem~1.1]{MN}. With Theorem~A and
Lemmas~\ref{lem:background} and~\ref{lem:WMN} this gives
$2\pi\le A(g)\le2\pi/\delta$, and $A(g)=2\pi/\delta$ only for the round metric of
curvature $\delta$.
\end{remark}
\section{The equality case without the injectivity radius hypothesis}\label{sec:conj}

In this section we discuss Conjecture~E. Proposition~\ref{prop:equiv} reformulates
it in terms of the injectivity radius. Lemma~\ref{lem:barrier} describes the
filling areas of loops shorter than $2\pi$. Corollary~\ref{cor:warped}
and Theorem~\ref{thm:leftinv} verify the conjecture for a class of doubly warped
metrics and for some left-invariant metrics. The lower bounds come from maps to the
round sphere (Lemma~\ref{lem:contract}).

The next proposition identifies the case of Conjecture~E that Theorem~B leaves
open.

\begin{proposition}\label{prop:equiv}
The following are equivalent.
\begin{enumerate}[label=\textup{(\alph*)}]
\item Conjecture~E.
\item Every metric on $S^3$ with $K_g\le1$ and a closed geodesic shorter than
$2\pi$ has $A(g)>2\pi$.
\item Every metric on $S^3$ with $K_g\le1$ and $A(g)=2\pi$ has $\inj(g)\ge\pi$.
\end{enumerate}
\end{proposition}
\begin{proof}
Assume (c). If $g$ has a closed geodesic shorter than $2\pi$, then $\inj(g)<\pi$
by Lemma~\ref{lem:kling}(a). So $A(g)\ne2\pi$ by (c), and $A(g)\ge2\pi$ by
Theorem~A. This is (b). Assume (b), and let $g$ be not round of any curvature,
$K=\max K_g>0$ and $g'=Kg$. Then $K_{g'}\le1$ and $A(g')=KA(g)\ge2\pi$
(Corollary~\ref{cor:round}). If all closed geodesics of $g'$ have length at
least $2\pi$, then $\inj(g')\ge\pi$ by Lemma~\ref{lem:kling}(a), and
$A(g')=2\pi$ would make $g'$ round by Theorem~B, and $g$ round of curvature $K$.
Otherwise (b) excludes it. So
$A(g)\max K_g>2\pi$. Assume (a). If $K_g\le1$ and $A(g)=2\pi$, then
$A(g)\max K_g\le2\pi$, so $g$ is round of some curvature $K\le1$ with
$2\pi/K=2\pi$; so $g$ is round and $\inj(g)=\pi$.
\end{proof}

If $K_g\le1$ and $\inj(g)<\pi$, then by Lemma~\ref{lem:kling}(b) a shortest
closed geodesic has length $2\inj(g)<2\pi$. It is embedded: otherwise it would split at a crossing $y$ into
two geodesic loops at $y$, each of length at least $2\inj(y)\ge2\inj(g)$. By
Theorem~\ref{thm:sharp}(c), a metric with $K_g\le1$ and $A(g)=2\pi$ that is not
round therefore has an embedded closed geodesic shorter than $2\pi$.

The next lemma describes the filling areas of loops shorter than $2\pi$.

\begin{lemma}\label{lem:barrier}
Let $(M,g)$ be a closed manifold of dimension at least three with $K_g\le1$,
and $c\in\Lambda M$ with $L(c)<2\pi$. Then
\[
 A(c,g)\le2\pi-\sqrt{4\pi^2-L(c)^2}\qquad\text{or}\qquad
 A(c,g)\ge2\pi+\sqrt{4\pi^2-L(c)^2}.
\]
A contractible embedded closed geodesic of length $\ell<2\pi$ has
$A\ge2\pi+\sqrt{4\pi^2-\ell^2}$. On a continuous family of loops in $\Lambda M$
of length at most $L_0<2\pi$, parametrized by a connected space and containing
a constant loop, the first alternative holds throughout.
\end{lemma}
\begin{proof}
A possibly branched minimal disk of area $a$ and boundary length $L$ in a
manifold with $K\le1$ satisfies $L^2\ge4\pi a-a^2$ \cite[Corollary~2.6]{MHH},
that is, $(a-2\pi)^2\ge4\pi^2-L^2$. Apply this to a least area disk spanning
a smooth embedded loop, which is smooth up to the boundary as in
Lemma~\ref{lem:F}(c). Smooth embedded loops are dense in $C^1$ in dimension at
least three, and length and filling area are continuous, so the inequality
holds for all loops. For an embedded closed geodesic Lemma~\ref{lem:GB} gives
$a\ge2\pi$. On a family as in the last statement $A$ is continuous, omits the
interval of radius $\sqrt{4\pi^2-L_0^2}$ about $2\pi$ and vanishes somewhere.
\end{proof}

\begin{remark}
Lemma~\ref{lem:barrier} rests on the inequality $L^2\ge4\pi a-a^2$ of
\cite{MHH} for minimal disks in manifolds with $K\le1$. Under the hypotheses of
Lemma~\ref{lem:barrier}, a loop of length at most $L_0<2\pi$ has
$|A(c,g)-2\pi|\ge\sqrt{4\pi^2-L_0^2}$. So on $S^3$, in a family of degree one
with $A(\Gamma,g)$ close to $2\pi$, the loops of almost maximal filling area are
long (Corollary~\ref{cor:barrier}).
\end{remark}

Together with Theorem~A, the lemma gives the following conclusions for families of
degree one.

\begin{corollary}\label{cor:barrier}
Let $g$ be a metric on $S^3$ with $K_g\le1$.
\begin{enumerate}[label=\textup{(\alph*)}]
\item Every family of degree one contains a loop of length at least $2\pi$.
\item If $A(\Gamma,g)\le2\pi+b$, $0\le b\le2\pi$, for a family of degree one, then
its loops
with $A\ge2\pi-b$ have length at least $\sqrt{4\pi^2-b^2}$.
\item Every family of degree one contains a loop with both length and filling
area at least $2\pi$.
\end{enumerate}
\end{corollary}
\begin{proof}
(a) Otherwise Lemma~\ref{lem:barrier} gives $A<2\pi$ on the family, against
Theorem~A. Part (b) follows from Lemma~\ref{lem:barrier}. (c) Otherwise, by the lemma, the closed set
$\{A\ge2\pi\}$ of the parameter disk equals the open set $\{L<2\pi,\ A>2\pi\}$;
it misses the boundary, so it is empty, against Theorem~A.
\end{proof}

\begin{remark}\label{rem:cylinder}
Mazet describes \cite[Section~1]{Mazet} a $C^{1,1}$ metric $\bar g$ on $S^3$. It is
the cylinder $S^2\times[-1,1]$ capped by two round hemispheres of the unit
three-sphere, with sectional curvature at most one in a weak sense and
Simon--Smith width $4\pi$. Mazet also describes smooth metrics $g_n\to\bar g$ with
$K_{g_n}\le1$ and $\widetilde W(g_n)\to4\pi$. For any metrics $g_n$ with
$K_{g_n}\le1$ and $\widetilde W(g_n)\to4\pi$, Corollary~\ref{cor:round} and
Lemmas~\ref{lem:background} and~\ref{lem:WMN} give
\[
 2\pi\le A(g_n)\max K_{g_n}\le A(g_n)\le\tfrac12W(g_n)\le
 \tfrac12\widetilde W(g_n)\longrightarrow2\pi .
\]
So for Mazet's metrics $A(g_n)\max K_{g_n}\to2\pi$, while the metrics converge to
$\bar g$, which is not round.
\end{remark}

For a linear map $L$ between Euclidean spaces, $|\Lambda^2L|$ is the norm of the
induced map on two-vectors, the product of the two largest singular values of
$L$ if there are two.

The next lemma gives lower bounds for widths from maps with bounded area dilation.
In particular a map of nonzero degree to the round sphere with area dilation at
most $B$ gives $A(g)\ge2\pi/B$ and $W(g)\ge4\pi/B$.

\begin{lemma}\label{lem:contract}
Let $G\colon(M,g)\to(N,h)$ be a $C^1$ map between closed Riemannian manifolds
with $|\Lambda^2dG|\le B$.
\begin{enumerate}[label=\textup{(\alph*)}]
\item For every $\alpha\in\pi_k(\Lambda M,M)$, $A(G_*\alpha,h)\le B\,A(\alpha,g)$.
\item If $M=N=S^3$, $h=g_0$ and $\deg G\ne0$, then $A(g)\ge2\pi/B$ and
$W(g)\ge4\pi/B$.
\end{enumerate}
\end{lemma}
\begin{proof}
(a) If $\Gamma$ represents $\alpha$, then $G\circ\Gamma$ represents $G_*\alpha$,
and if a Lipschitz disk $f$ spans a loop $c$, then $G\circ f$ spans $G\circ c$
and has $h$-area at most $B\Area_g(f)$. (b) If $\alpha$ has degree one,
$G_*\alpha$ has degree $\pm\deg G\ne0$, so it is nontrivial, and
$A(G_*\alpha,g_0)\ge2\pi$ by Theorem~A. Apply (a) and $2A\le W$.
\end{proof}

\begin{remark}
By Lemma~\ref{lem:contract}(b), a map of nonzero degree onto the round sphere
with $|\Lambda^2dG|\le B<1$ gives $A(g)>2\pi$. The map may increase lengths.
\end{remark}

A \emph{doubly warped metric} on $S^3$ is a metric
$g=dt^2+a(t)^2dx^2+b(t)^2dy^2$, $0\le t\le T$, with $x,y$ of period $2\pi$,
$b(0)=0$, $b'(0)=1$, $a(T)=0$, $a'(T)=-1$, $a>0$ on $[0,T)$ and $b>0$ on
$(0,T]$, that extends smoothly over the cores $\{t=0\}$ and $\{t=T\}$. Then also
$a'(0)=b'(T)=0$. Its curvature operator is diagonal in the coordinate frame, with
eigenvalues $-a''/a$, $-b''/b$ and $-a'b'/(ab)$, so $K_g\le1$ means
$a''+a\ge0$, $b''+b\ge0$ and $a'b'+ab\ge0$. For $T=\frac\pi2$, $a=\cos t$,
$b=\sin t$ it is the round metric, in the coordinates
$(t,x,y)\mapsto(\cos t\,e^{ix},\sin t\,e^{iy})$ of $S^3\subset\C^2$.

The next theorem gives a lower bound for doubly warped metrics with $T\le\pi$ that
satisfy the two radial inequalities $a''+a\ge0$ and $b''+b\ge0$.

\begin{theorem}\label{thm:warped}
Let $g=dt^2+a^2dx^2+b^2dy^2$ be a doubly warped metric on $S^3$ with
$a''+a\ge0$, $b''+b\ge0$ and $T\le\pi$. Then $T\ge\frac\pi2$, with equality
only for the round metric, and
\[
 A(g)\ \ge\ 2\sqrt2\,\pi\,\bigl(\cos\tfrac T2-\cos T\bigr).
\]
In particular $A(g)>2\pi$ unless $g$ is round.
\end{theorem}
\begin{proof}
The function $h=b'\sin t-b\cos t$ has $h(0)=0$ and
$h'=(b''+b)\sin t\ge0$ on $[0,T]$, so $b/\sin t$ is nondecreasing with limit
$1$ at zero, and $b\ge\sin t$; likewise $a\ge\sin(T-t)$. Since $T\le\pi$ and
$b'(T)=0$, $h(T)=-b(T)\cos T\ge0$, so $T\ge\frac\pi2$; if $T=\frac\pi2$, then
$h=0$, $b=\sin t$, and likewise $a=\cos t$.

For $\frac\pi2<T\le\pi$ put $k=-\cos T$, $P(t)=\cos(T-t)-\cos T$
and $c=\sqrt2(\cos\frac T2-\cos T)$. Define $f\colon[0,T]\to[0,\frac\pi2]$
by $c\sin f=P$ on $[0,\frac T2]$ and $f(t)=\frac\pi2-f(T-t)$ on
$[\frac T2,T]$; then $f(\frac T2)=\frac\pi4$, and $f$ is $C^1$ and
increasing. The map $(t,x,y)\mapsto(f(t),x,y)$ onto the round
sphere $df^2+\cos^2f\,dx^2+\sin^2f\,dy^2$ is a homeomorphism,
$C^1$ also at the cores. Its singular values are $f'$, $\cos f/a$ and
$\sin f/b$. They are at most $f'$, $\cos f/a_0$ and $\sin f/b_0$, where
$a_0=\sin(T-t)\le a$ and $b_0=\sin t\le b$. The
three products of two of them are at most $1/c$ on $[0,\frac T2]$, hence, by
the symmetry $t\mapsto T-t$, on $[0,T]$. The quotients below are taken in the
interior, and the bounds extend to the endpoints by continuity of the
differential.

First, $f'\cos f=(\sin f)'=a_0/c$, so
$f'\cos f/a\le1/c$.

Second, $f'\sin f/b_0=a_0\tan f/(cb_0)$, and
$\tan f\le b_0/a_0$ is equivalent to $H\le c$, where
$H=P\sqrt{a_0^2+b_0^2}/b_0$. With
$v=\frac t2$ one has $P/b_0=\sin(T-v)/\cos v$ and
$a_0^2+b_0^2=1+k\cos(T-4v)$, so
\[
 \tfrac12(\log H^2)'=\tan v-\cot(T-v)+\frac{2k\sin(T-4v)}{1+k\cos(T-4v)}
 \ \ge0\qquad\text{on }(0,\tfrac T4],
\]
because $T\ge\frac\pi2$ and $0\le T-4v<\pi$. Hence $H\le H(\frac T2)=c$.

Third, $\sin f\cos f/(a_0b_0)\le1/c$ means
$c\sin2f\le2a_0b_0=P(2t)$. With $s=T-\frac\pi2$,
$P(z)=k+\sin(z-s)$ is increasing on $[0,T]$, and
$t=s+\arcsin(c\sin f-k)$. So the inequality reads $D(f)\ge0$ on
$[0,\frac\pi4]$, where
\[
 D(f)=s+2\arcsin(c\sin f-k)-\arcsin(c\sin2f-k),
\]
and $D(0)=s-\arcsin k=0$. The arguments of $\arcsin$ lie in
$[-k,c-k]\subset[-1,1]$, since $1-c+k=(2\sqrt2-2)(\frac1{\sqrt2}-\cos\frac
T2)(1-\cos\frac T2)\ge0$. The inequality $D'\ge0$ is equivalent to
\[
 \cos^2f\,\bigl(1-(c\sin2f-k)^2\bigr)
 -\cos^22f\,\bigl(1-(c\sin f-k)^2\bigr)=\sin f\cdot Q\ \ge0,
\]
\[
 Q=2ck\,E-(c^2+k^2-1)\sin3f,\qquad E=2\cos^3f-\cos^22f .
\]
With $m=\cos\frac T2\in[0,\frac1{\sqrt2})$ one has $k=1-2m^2$,
$c=\sqrt2(1-m)(1+2m)$ and $c^2+k^2-1=2(1-m)(1+3m)k$, so
\[
 Q=2(1-m)k\,\bigl[\sqrt2(1+2m)E-(1+3m)\sin3f\bigr].
\]
The bracket is affine in $m$, and at $m=0$ and $m=\frac1{\sqrt2}$ it is
nonnegative if $E\ge(\sqrt2-\frac12)\sin3f$. Put
$e=\sqrt2\cos f\in[1,\sqrt2]$. Then $E=e^3/\sqrt2-(e^2-1)^2$ and
$\sin3f=\sqrt{2-e^2}\,(2e^2-1)/\sqrt2$, and with
$\sqrt{2-e^2}\le\frac12(3-e^2)$,
\[
 4\sqrt2\,\bigl(E-(\sqrt2-\tfrac12)\sin3f\bigr)\ \ge\ -2e^4+4e^3-(6\sqrt2-7)e^2-(3-2\sqrt2),
\]
which is concave on $[1,\sqrt2]$ with the values $6-4\sqrt2>0$ and
$3-2\sqrt2>0$ at the ends.

So $|\Lambda^2|\le1/c$, and Lemma~\ref{lem:contract}(b) gives
$A(g)\ge2\pi c$. Finally $c=\sqrt2(1+m-2m^2)$ is concave in $m$ with the values
$\sqrt2$ at $m=0$ and $1$ at $m=\frac1{\sqrt2}$, so $c>1$ for $T>\frac\pi2$.
\end{proof}

Together with a scaling, the theorem proves Conjecture~E for these metrics.

\begin{corollary}\label{cor:warped}
Conjecture~E holds for the metrics of Theorem~\ref{thm:warped}: they satisfy
$A(g)\max K_g>2\pi$ unless $g$ is round of some curvature.
\end{corollary}
\begin{proof}
If $\max K_g\ge1$, then $A(g)\max K_g\ge A(g)>2\pi$ unless $g$ is round. If
$\max K_g<1$, the metric $(\max K_g)\,g$ is doubly warped, with curvature at
most one and length $T\sqrt{\max K_g}\le\pi$, and Theorem~\ref{thm:warped}
applies to it.
\end{proof}

The next theorem treats left-invariant metrics on $SU(2)$.

\begin{theorem}\label{thm:leftinv}
Let $X_1,X_2,X_3$ be a left-invariant frame on $S^3=SU(2)$, orthonormal for
$g_0$, with dual coframe $s_1,s_2,s_3$, and $g=a_1^2s_1^2+a_2^2s_2^2+a_3^2s_3^2$
with $0<a_1\le a_2\le a_3$. Then $A(g)\ge2\pi a_1a_2$. If moreover
$a_1^2\ge\frac14a_3^2$, then $A(g)\max K_g>2\pi$ unless $a_1=a_2=a_3$.
\end{theorem}
\begin{proof}
The identity map from $g$ to $g_0$ has singular values $1/a_i$, so
$|\Lambda^2|\le1/(a_1a_2)$, and Lemma~\ref{lem:contract}(b) gives the first
bound. Put $x_i=a_i^2$. By Milnor's formulas \cite[\S4]{Milnor} the curvature
operator of $g$ is diagonal in the frame $X_i\wedge X_j$, with
\[
 K_{ij}=\frac{(x_i-x_j)^2+x_k(2x_i+2x_j-3x_k)}{x_1x_2x_3},
 \qquad
 K_{ij}-K_{ik}=\frac{4(x_k-x_j)(x_i-x_j-x_k)}{x_1x_2x_3}
\]
for $\{i,j,k\}=\{1,2,3\}$. So $K_{12}$ is the smallest, and the largest is
$K_{23}$ if $x_3\le x_1+x_2$ and $K_{13}$ otherwise. The product $A\max K$ is
scale invariant; scale to $x_3=1$ and write $p=x_1\le q=x_2\le1$. Then
$\max K=N/(pq)$ with $N=(1-q)^2+p(2q+2-3p)$ if $p+q\ge1$ and
$N=(1-p)^2+q(2p+2-3q)$ if $p+q<1$. In the second case $N$ is concave in
$q\in[p,1-p]$ with end values $1$ and $4p(1-p)$, so $N\ge4p(1-p)$.

If $p\ge\frac14$, the first bound gives $A\max K\ge2\pi N/\sqrt{pq}$, and
$N\ge\sqrt{pq}$ with equality only for $p=q=1$. Indeed, for $p+q\ge1$,
$N\ge p(2q+2-3p)$, and $\sqrt p\,(2q+2-3p)-\sqrt q$ is nondecreasing in
$q\ge p$ with the value $\sqrt p\,(1-p)\ge0$ at $q=p$. For $p+q<1$,
$N\ge4p(1-p)>\sqrt{p(1-p)}>\sqrt{pq}$ since $p(1-p)\ge\frac3{16}$.
\end{proof}

\begin{remark}
Theorem~\ref{thm:leftinv} applies to metrics not covered by Theorem~B. Let
$a_2=a_3=1$ and $\frac14\le a_1^2<\frac13$. By the formulas above,
$\max K_g=4-3a_1^2$. The integral curves of $X_1$ are closed geodesics of $g$ of
length $2\pi a_1$, since $[X_1,Y]$ is $g$-orthogonal to $X_1$ for every left-invariant $Y$. In the
metric $(\max K_g)\,g$ they have length $2\pi a_1\sqrt{4-3a_1^2}<2\pi$, so its
injectivity radius is less than $\pi$ by Lemma~\ref{lem:kling}.
\end{remark}

\end{document}